# A WAITING TIME PROBLEM ARISING FROM THE STUDY OF MULTI-STAGE CARCINOGENESIS

By Rick Durrett,[1] Deena Schmidt[1,2] and Jason Schweinsberg[3]

*Cornell University, University of Minnesota and
University of California, San Diego*

We consider the population genetics problem: how long does it take before some member of the population has $m$ specified mutations? The case $m = 2$ is relevant to onset of cancer due to the inactivation of both copies of a tumor suppressor gene. Models for larger $m$ are needed for colon cancer and other diseases where a sequence of mutations leads to cells with uncontrolled growth.

**1. Introduction.** It has long been known that cancer is the end result of several mutations that disrupt normal cell division. Armitage and Doll [1] did a statistical analysis of the age of onset of several cancers and fit power laws to estimate the number of mutations. Knudson [15] discovered that the incidence of retinoblastoma (cancer of the retina) grows as a linear function of time in the group of children who have multiple cancers in both eyes, but as a slower quadratic function in children who only have one cancer. Based on this, Knudson proposed the concept of a tumor suppressor gene. Later it was confirmed that in the first group of children, one copy is already inactivated at birth, while in the second group both copies must be mutated before cancer occurs. Since that time, about 30 tumor suppressor genes have been identified. They have the property that inactivating the first copy does not cause a change, while inactivating the second increases the cells' net reproductive rate, which is a step toward cancer.

There is now considerable evidence that colon cancer is the end result of several mutations. The earliest evidence was statistical. Luebeck and Mool-

Received July 2007; revised July 2008.
[1]Supported in part by NSF Grant DMS-02-02935 from the probability program and NSF/NIGMS Grant DMS-02-01037.
[2]Supported by an NSF graduate fellow.
[3]Supported in part by NSF Grants DMS-05-04882 and DMS-08-05472.

gavakar [18] fit a four-stage model to the age-specific incidence of colorectal cancers in the Surveillance, Epidemiology, and End Results registry, which cover 10 percent of the US population. Calabrese et al. [5] examined 1022 colorectal cancers sampled from nine large regional hospitals in southeastern Finland. They found support for a model with five or six oncogenic mutations in individuals with hereditary risk factors and seven or eight mutations in patients without.

Over the last decade, a number of studies have been carried out to identify the molecular pathways involved in the development of colorectal cancer. See Jones et al. [14] for a recent report. The process is initiated when a single colorectal cell acquires mutations inactivating the ACP/$\beta$-catenin pathway. This results in the growth of small benign tumor (adenoma). Subsequent mutations in a short list of other pathways transform the adenoma into a malignant tumor (carcinoma), and lead to metastasis, the ability of the cancer to spread to other organs.

In this paper, we propose a simple mathematical model for cancer development in which cancer occurs when one cell accumulates $m$ mutations. Consider a population of fixed size $N$. Readers who are used to the study of the genetics of diploid organisms may have expected to see $2N$ here, but our concern is for a collection of $N$ cells. We choose a model in which the number of cells is fixed because organs in the body are typically of constant size. We assume that the population evolves according to the Moran model, which was first proposed by Moran [19]. That is, each individual lives for an exponentially distributed amount of time with mean one, and then is replaced by a new individual whose parent is chosen at random from the $N$ individuals in the population (including the one being replaced). For more on this model, see Section 3.4 of [11].

In our model, each individual has a type $0 \leq j \leq m$. Initially, all individuals have type 0. In the usual population genetics model, mutations only occur at replacement events. We assume instead that types are clonally inherited, that is, every individual has the same type as its parent. However, thinking of a collection of cells that may acquire mutations due to radiation or other environmental factors, we will suppose that during their lifetimes, individuals of type $j-1$ mutate to type $j$ at rate $u_j$. We call such a mutation a *type $j$ mutation*. Let $X_j(t)$ be the number of type $j$ individuals at time $t$. For each positive integer $m$, let $\tau_m = \inf\{t : X_m(t) > 0\}$ be the first time at which there is an individual in the population of type $m$. Clearly, $\tau_1$ has the exponential distribution with rate $Nu_1$. Our goal is to compute the asymptotic distribution of $\tau_m$ for $m \geq 2$ as $N \to \infty$.

We begin by considering the case $m = 2$ and discussing previous work. Schinazi [21, 22] has considered related questions. In the first paper, he computes the probability that in a branching process where individuals have two offspring with probability $p$ and zero with probability $1-p$, a mutation will arise before the process dies out. In the second paper, he uses this to investigate the probability of a type 2 mutation when type 0 cells divide a fixed number of times with the possibility of mutating to a type 1 cell that begins a binary branching process.

More relevant to our investigation is the work of Komarova, Sengupta and Nowak [17], Iwasa, Michor and Nowak [13] and Iwasa et al. [12]. Their analysis begins with the observation that while the number of mutant individuals is $o(N)$, we can approximate the number of cells with mutations by a branching process in which each individual gives birth at rate one and dies at rate one. Let $Z$ denote the total progeny of such a branching process. Since the embedded discrete time Markov chain is a simple random walk, we have (see, e.g., page 197 in [7])

$$P(Z > n) = 2^{-2n} \binom{2n}{n} \sim \frac{1}{\sqrt{\pi n}}.$$

If we ignore interference between successive new type 1 mutations, then their total progeny $Z_1, Z_2, \ldots$ are i.i.d. variables in the domain of attraction of a stable law with index $1/2$, so $\max_{i \leq M} Z_i$ and $Z_1 + \cdots + Z_M$ will be $O(M^2)$. Therefore, we expect to see our first type 2 mutation in the family of the $M$th type 1 mutation, where $M = O(1/\sqrt{u_2})$. Standard results for simple random walk imply that the largest of our first $M$ families will have $O(M)$ type 1 individuals alive at the same time, so for the branching process approximation to hold, we need $1/\sqrt{u_2} \ll N$, where here and throughout the paper, $f(N) \ll g(N)$ means that $f(N)/g(N) \to 0$ as $N \to \infty$. Type 1 mutations occur at rate $Nu_1$, so a type 2 mutation will first occur at a time of order $1/Nu_1\sqrt{u_2}$.

As long as the branching process approximation is accurate, the amount of time we have to wait for a type 1 mutation that will have a type 2 individual as a descendant will be approximately exponential, since mutations occur at times of a Poisson process with rate $Nu_1$ and the type 1 mutations that lead to a type 2 are a thinning of that process in which points are kept with probability $\sim \sqrt{u_2}$, which is $O(1/M)$, where here and throughout the paper, $f(N) \sim g(N)$ means that $f(N)/g(N) \to 1$ as $N \to \infty$. The duration of the longest of $M$ type 1 families is $O(M)$, so the time between when the type 1 mutation occurs and when the type 2 descendant appears is $O(1/\sqrt{u_2})$. This will be negligible in comparison to $1/Nu_1\sqrt{u_2}$ as long as $Nu_1 \ll 1$, so the waiting time for the first type 2 individual will also be approximately exponential. This leads to a result stated on pages 231–232 of Nowak's book [20] on Evolutionary Dynamics. If $1/\sqrt{u_2} \ll N \ll 1/u_1$, then

$$(1.1) \qquad P(\tau_2 \leq t) \approx 1 - \exp(-Nu_1\sqrt{u_2}t).$$

Figure 1 shows the distribution of $\tau_2 \cdot Nu_1\sqrt{u_2}$ in 10,000 simulations of the Moran model when $N = 10^3$ and $u_1 = u_2 = 10^{-4}$. Here, $Nu_1 = 0.1$ and $N\sqrt{u_2} = 10$, so as the last result predicts, the scaled waiting time is approximately exponential.

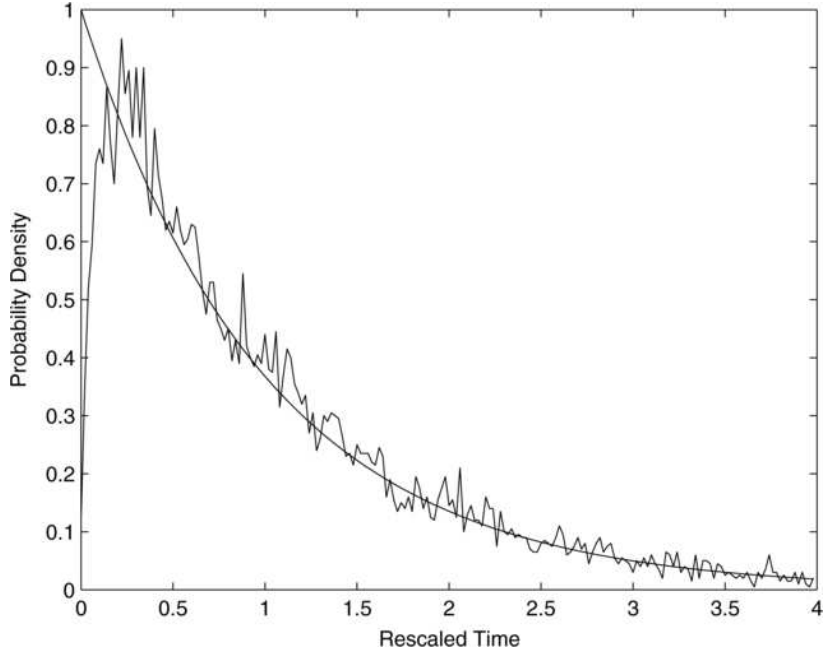

Fig. 1. *Distribution of $\tau_2 \cdot Nu_1\sqrt{u_2} = 1000$ in 10,000 simulations when $N = 10^3$ and $u_1 = u_2 = 10^{-4}$. $Nu_1 = 0.1$ and $N\sqrt{u_2} = 10$, so as (1.1) predicts the scaled waiting time is approximately exponential.*

We do not refer to the result given in (1.1) as a theorem because their argument is not completely rigorous. For example, the authors use the branching process approximation without proving it is valid. However, this is a minor quibble, since as the reader will see in Section 2, it is straightforward to fill in the missing details and establish the following more general result.

THEOREM 1. *Suppose that $Nu_1 \to \lambda \in [0, \infty)$, $u_2 \to 0$ and $N\sqrt{u_2} \to \infty$ as $N \to \infty$. Then $\tau_2 \cdot Nu_1\sqrt{u_2}$ converges to a limit that has density function*

$$f_2(t) = h(t) \exp\left(-\int_0^t h(s)\,ds\right) \qquad \text{where } h(s) = \frac{1 - e^{-2s/\lambda}}{1 + e^{-2s/\lambda}},$$

*if $\lambda > 0$ and $f_2(t) = e^{-t}$ if $\lambda = 0$.*

Here, $h(t)$ is the hazard function, that is, if we let $F_2(t) = \exp(-\int_0^t h(s)\,ds)$ be the tail of the distribution, then $h(t) = f_2(t)/F_2(t)$. Figure 2 shows the distribution of $\tau_2 \cdot Nu_1\sqrt{u_2}$ in 10,000 simulations of the Moran model when $N = 10^3$, $u_1 = 10^{-3}$ and $u_2 = 10^{-4}$. $Nu_1 = 1$ so the limit is not exponential, but Theorem 1 gives a good fit to the observed distribution.

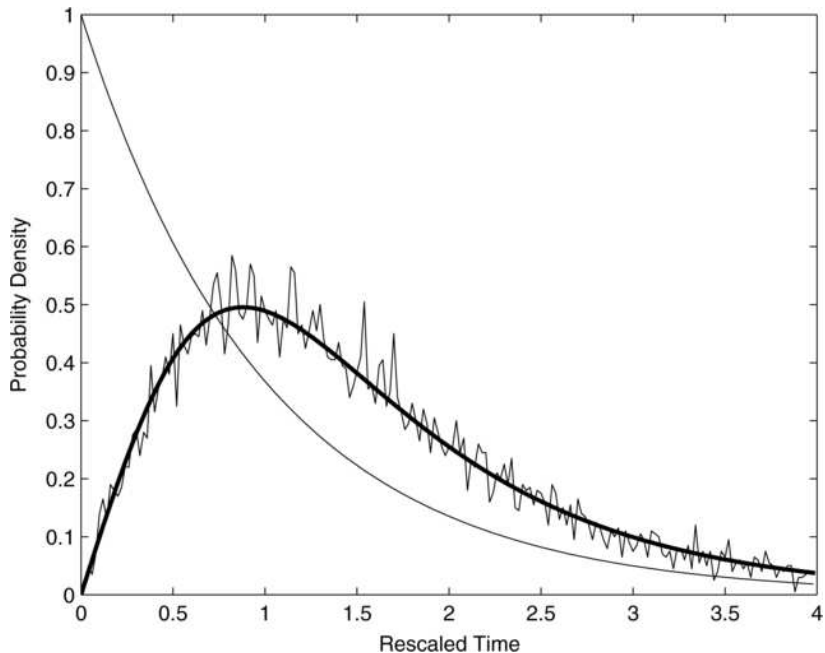

Fig. 2. *Distribution of $\tau_2 \cdot Nu_1\sqrt{u_2} = 1000$ in 10,000 simulations when $N = 10^3$, $u_1 = 10^{-3}$ and $u_2 = 10^{-4}$. $Nu_1 = 1$ and $N\sqrt{u_2} = 0.1$, so the limit is not exponential, but is fit well by the result in Theorem 1.*

Before turning to the case of $m$ mutations, we should clarify one point. In our model, mutations occur during the lifetime of an individual, but in the following discussion, we will count births to estimate the probability a desired mutation will occur. This might seem to only be appropriate if mutations occur at birth. However, since each individual lives for an exponential amount of time with mean 1, the number of "man-hours" $\int_0^{T_0} X_1(s)\,ds$ before the family dies out at time $T_0$ is roughly the same as the number of births. In any case, the following discussion is only a heuristic that helps explain the answer, but does not directly enter into its proof.

To extend the analysis to the $m$-stage waiting time problem, suppose $M$ distinct type 1 mutations have appeared. If the family sizes of these $M$ mutations can be modeled by independent branching processes, the total number of offspring of type 1 individuals will be $O(M^2)$. Because each type 1 individual mutates to type 2 at rate $u_2$, there will be $O(M^2 u_2)$ mutations that produce type 2 individuals. The total progeny of these individuals will consist of $O(M^4 u_2^2)$ type 2 individuals. We can expect to see our first type 3 individual when $M^4 u_2^2 = O(1/u_3)$ or $M = O(u_2^{-1/2} u_3^{-1/4})$. Thus, for the branching process approximation to hold, we need $u_2^{-1/2} u_3^{-1/4} \ll N$. Since type 1 mutations occur at rate $Nu_1$, the expected waiting time will be of order

$$1/Nu_1 u_2^{1/2} u_3^{1/4}.$$

To help develop a good mental picture, it is instructive to consider the numerical example in which $N = 10^5$, $u_1 = 10^{-6}$, $u_2 = 10^{-5}$ and $u_3 = 10^{-4}$. By the reasoning above, we will first see a type 3 mutation when the number of type 2's is of order $100 = 1/\sqrt{u_3}$, since in this case there will be of order $10{,}000 = 1/u_3$ type 2 births before the family dies out. To have a type 2 family reach size 100, we will need 100 mutations from type 1 to type 2, and for this we will need of order $100/u_2 = 10^7$ type 1 births, which will in turn occur if the type 1 family reaches size of order $10^{7/2} \approx 3162$. Note that $X_2(t) \ll X_1(t)$ and within the time that the large type 1 family exists, 100's of type 2 families will be started and die out. This difference in the time and size scales for the processes $X_i(t)$ is a complicating factor in the proof, but ultimately it also allows us to separate the type 1's from types 2 to $m$ and use induction.

Extrapolating the calculation above to $m$ stages, we let

$$(1.2) \qquad r_{j,m} = u_{j+1}^{1/2} u_{j+2}^{1/4} \cdots u_m^{1/2^{m-j}}$$

for $1 \leq j < m$, and set $r_{m,m} = 1$ and $r_{0,m} = u_1 r_{1,m}$. Let $q_{j,m}$ be the probability a type $j$ individual gives rise to a type $m$ descendant. We will show that $q_{j,m} \sim r_{j,m}$, so we will need of order $1/r_{j,m}$ mutations to type $j$ before time $\tau_m$.

THEOREM 2. *Fix an integer $m \geq 2$. Suppose that:*

(i) $Nu_1 \to 0$.
(ii) *For $j = 1, \ldots, m-1$, there is a constant $b_j > 0$ such that $u_{j+1}/u_j > b_j$ for all $N$.*
(iii) *There is an $a > 0$ so that $N^a u_m \to 0$.*
(iv) $Nr_{1,m} \to \infty$.
*Then for all $t > 0$,*

$$(1.3) \qquad \lim_{N \to \infty} P(\tau_m > t/Nr_{0,m}) = \exp(-t).$$

As discussed above, condition (iv) which says $1/r_{1,m} \ll N$ is needed for the branching process assumption to be valid, and condition (i) is needed for the waiting time to be exponential, because if (i) fails then the time between the type 1 mutation that will have a type $m$ descendant and the birth of the type $m$ descendant cannot be neglected. If $u_j = \mu$ for all $j$, (ii) is trivial. In this case $r_{1,m} = \mu^{a(m)}$, where $a(m) = 1 - 2^{-(m-1)}$. Conditions (i) and (iv) become $N^{-1/a(m)} \ll \mu \ll N^{-1}$, and when condition (i) is satisfied, (iii) holds.

Conditions (ii) and (iii) are technicalities that allow us to prove the result without having to suppose that $u_j \equiv \mu$, which would not be natural in modeling cancer. In the presence of (ii), condition (iii) ensures that $\max_{j \leq m} u_j \ll N^{-a}$ for some $a > 0$. This is natural because even in the late stages of progression to cancer, the per cell division mutation probabilities are small.

Condition (ii) is motivated by the fact that in most cancers we expect $u_j$ to be increasing in $j$. The simple extension of this given in (ii) is useful so that we do not rule out some interesting special cases. In modeling the tumor suppressor genes mentioned earlier, it is natural to take $u_1 = 2\mu$ and $u_2 = \mu$, that is, at the first stage a mutation can knock out one of the two copies of the gene, but after this occurs, there is only one copy subject to mutation. A case with $u_1/u_2 = 30$ occurs in Durrett and Schmidt's study of regulatory sequence evolution [9].

Condition (iv) ensures that an individual of type $m$ will appear before any type 1 mutation achieves fixation. In the case $m = 2$, Iwasa et al. [13] called this *stochastic tunneling*. A given type 1 mutation fixates with probability $1/N$ and type 1 mutations occur at rate approximately $Nu_1$, so fixation occurs before a type $m$ individual appears if $Nr_{1,m} \to 0$, and then once a type 1 mutation fixates, the problem reduces to the problem of waiting for $m-1$ additional mutations. In the borderline case considered in the next result, either a type $m$ individual could appear before fixation, or a type $m$ mutation could be achieved through the fixation of type 1 individuals followed by the generation of an individual with $m-1$ additional mutations.

THEOREM 3. *Fix an integer $m \geq 2$. Assume conditions* (i), (ii) *and* (iii) *from Theorem 2 hold. If $(Nr_{1,m})^2 \to \gamma > 0$, and we let*

$$(1.4) \qquad \alpha = \sum_{k=1}^{\infty} \frac{\gamma^k}{(k-1)!(k-1)!} \bigg/ \sum_{k=1}^{\infty} \frac{\gamma^k}{k!(k-1)!} > 1,$$

*then for all $t > 0$, $\lim_{N \to \infty} P(u_1 \tau_m > t) = \exp(-\alpha t)$.*

Figure 3 shows the distribution of $u_1 \tau_2$ in 10,000 simulations of the Moran model when $N = 10^3$, $u_1 = 10^{-4}$ and $u_2 = 10^{-6}$. $Nu_1 = 0.1$ and $N\sqrt{u_2} = 1$, so the assumptions of Theorem 3 hold with $\gamma = 1$. Numerically evaluating the constant gives $\alpha = 1.433$ and as the figure shows the exponential with this rate gives a reasonable fit to the simulated data.

Theorem 3 will be proved by reducing the general case to a two-type model with $\bar{u}_1 = u_1$ and $\bar{u}_2 = u_2 q_{2,m} \sim r_{1,m}^2$. We will show that it suffices to do calculations for a model in which type 1 mutations are not allowed when the number of type 1 individuals $X_1(t)$ is positive. In this case, if we start with $X_1(0) = N\varepsilon$ then $N^{-1}X_1(Nt) \to Z_t$ where $Z_t$ is the Wright–Fisher diffusion process with infinitesimal generator $x(1-x)d^2/dx^2$. When $X_1(Nt) = Nx$, mutations to type 2 that eventually lead to a type $m$ individual occur at rate approximately

$$N \cdot Nx \cdot u_2 q_{2,m} \sim N^2 r_{1,m}^2 x \to \gamma x,$$

so, if we let $u(x)$ be the probability that the process $Z_t$ hits 0 before reaching 1 or generating a type $m$ mutation, then $u(x)$ satisfies

(1.5) $\qquad x(1-x)u''(x) - \gamma x u(x) = 0, \qquad u(0) = 1, \qquad u(1) = 0.$

The constant $\alpha = \lim_{\varepsilon \to 0}(1 - u(\varepsilon))/\varepsilon$. Its relevance for the problem is that starting from a single type 1 individual, the probability of reaching $N$ or generating a type $m$ mutation is $\sim \alpha/N$. Since mutations to type 1 occur at rate $\sim Nu_1$, the waiting time is roughly exponential with rate $u_1 \alpha$.

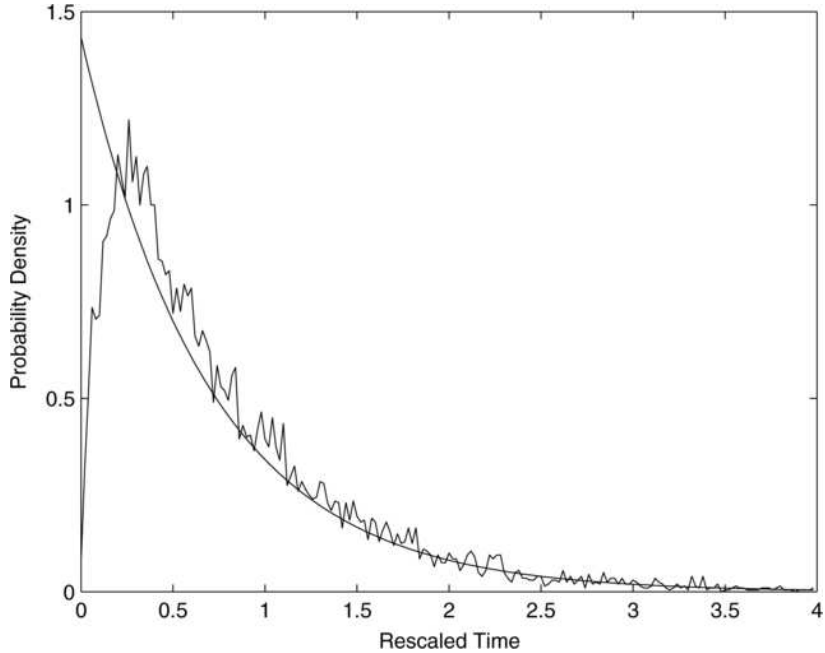

FIG. 3. *Distribution of $u_1\tau_2$ when $N = 10^3$, $u_1 = 10^{-4}$ and $u_2 = 10^{-6}$. $Nu_1 = 0.1$ and $N\sqrt{u_2} = 1$ so we are in the regime covered by Theorem 3. The constant $\gamma = 1$ so $\alpha = 1.433$. As the graph shows the exponential distribution with rate $\alpha$ gives a reasonbly good fit to the simulated data.*

One can check (see Lemma 6.9 below) that (1.5) can be solved by the following power series around $x = 1$:

$$(1.6) \quad u(x) = c \sum_{k=1}^{\infty} \frac{\gamma^k}{k!(k-1)!}(1-x)^k.$$

Picking $c$ so that $u(0) = 1$, it follows that $\alpha$ has the form given in (1.4). Another approach to solving (1.5) is to use the Feynman–Kac formula; see formula (3.19.5b) on page 225 of [4].

We do not discuss in this paper the case $Nu_1 \to \infty$. We instead refer the reader to [23], where asymptotic results in this regime are obtained in the special case when $u_j = \mu$ for all $j$.

The rest of this paper is organized as follows. In Section 2, we give the proof of Theorem 1. In Section 3, we collect some results for a two-type population model that will be useful later in the paper. In Section 4, we calculate by induction the probability that a given type 1 individual has a type $m$ descendant. In Section 5, we combine this result with a Poisson approximation result of Arratia, Goldstein and Gordon [2] to prove Theorem 2. Theorem 3 is proved in Sections 6 and 7. Throughout our proofs, $C$ denotes a constant whose value is unimportant and will change from line to line.

**2. Proof of Theorem 1.** If we let $X_1(t)$ be the number of type 1 individuals at time $t$ then

$$(2.1) \quad P(\tau_2 > t) = E\exp\left(-u_2 \int_0^t X_1(s)\,ds\right),$$

because at time $s$, there are $X_1(s)$ individuals each experiencing type 2 mutations at rate $u_2$. We will compare $X_1(t)$ with a continuous-time branching process with immigration, $Y(t)$. When $X_1(t) = k$, type 1 mutations occur at rate $(N-k)u_1$, while birth events in which a type 1 individual replaces a type 0 individual occur at rate $k(N-k)/N$, so before time $\tau_2$, we have jumps

$$k \to k+1 \quad \text{at rate } (k + Nu_1) \cdot \frac{N-k}{N},$$

$$k \to k-1 \quad \text{at rate } k \cdot \frac{N-k}{N}.$$

In the branching process with immigration, $Y(t)$, we have jumps

$$k \to k+1 \quad \text{at rate } k + Nu_1,$$

$$k \to k-1 \quad \text{at rate } k.$$

Therefore, up to time $\tau_2$, the process $\{X_1(t), t \geq 0\}$ is a time-change of $\{Y(t), t \geq 0\}$, in which time runs slower than in the branching process by a factor of $(N - X_1(t))/N$. That is, if

$$T(t) = \int_0^t \frac{N - X_1(s)}{N} \, ds \leq t,$$

then the two processes can be coupled so that $X_1(t) = Y(T(t))$, for all $t \geq 0$. The time change will have little effect as long as $X_1(t)$ is $o(N)$. The next lemma shows that on the relevant time scale, the number of mutants stays small with high probability.

LEMMA 2.1.  *Fix $t > 0$, $\varepsilon > 0$, and let $M_t = \max_{0 \leq s \leq t/(Nu_1\sqrt{u_2})} X_1(s)$. We have*

$$\lim_{N \to \infty} P(M_t > \varepsilon N) = 0.$$

PROOF.  Since mutant individuals give birth and die at the same rate, the process $\{X_1(s), s \geq 0\}$ is a submartingale. Because the rate of type 1 mutations is always bounded above by $Nu_1$, we have $EX_1(s) \leq Nu_1 s$ for all $s$. By Doob's maximal inequality,

$$P(M_t > \varepsilon N) \leq \frac{EX_1(t/Nu_1\sqrt{u_2})}{\varepsilon N} \leq \frac{t}{\varepsilon N \sqrt{u_2}},$$

which goes to zero as $N \to \infty$, since $N\sqrt{u_2} \to \infty$.  □

Using the time change in (2.1), we have

$$P(\tau_2 > t/Nu_1\sqrt{u_2}) = E \exp\left(-u_2 \int_0^{t/Nu_1\sqrt{u_2}} Y(T(s)) \, ds\right).$$

Changing variables $r = T(s)$, which means $s = U(r)$, where $U = T^{-1}$, $ds = U'(r) \, dr$ and the above is

$$P(\tau_2 > t/Nu_1\sqrt{u_2}) = E \exp\left(-u_2 \int_0^{T(t/Nu_1\sqrt{u_2})} Y(r) U'(r) \, dr\right).$$

When $M_t \leq N\varepsilon$, $1 \geq T'(t) \geq 1 - \varepsilon$, so the inverse function has slope $1 \leq U'(r) \leq 1/(1 - \varepsilon)$. Thus, in view of Lemma 2.1, it is enough to prove the result for the branching process, $Y(t)$.

Use $Q$ to denote the distribution of $\{Y(t), t \geq 0\}$, and let $Q_1$ denote the law of the process starting from a single type 1 and modified to have no further mutations to type 1. We first compute $g_2(t) = Q_1(\tau_2 \leq t)$. Wodarz and Komarova [24] do this, see pages 37–39, by using Kolmogorov's forward equation to get a partial differential equation

$$\frac{\partial \phi}{\partial t}(t, y) = (y^2 - (2 + u_2)y + 1)\frac{\partial \phi}{\partial y}(t, y)$$

for the generating function $\phi(t,y) = \sum_j Q_1(X_1(t) = j, X_2(t) = 0)y^j$ of the system in which type 2's are not allowed to give birth or die. They use the method of characteristics to reduce the PDE to a Riccati ordinary differential equation. To help readers who want to follow their derivation, we note that the last equation on page 38 is missing a factor of $j$ in the last term and in the change of variables from $y$ to $z$ on page 39, 2 should be $r$.

Here, we will use Kolmogorov's backward differential equation to derive an ODE, which has the advantage that it generalizes easily to the $m$ stage problem. By considering what happens between time 0 and $h$,

$$g_2(t+h) = g_2(t)[1 - (2+u_2)h] + h[2g_2(t) - g_2(t)^2] + h \cdot 0 + u_2 h \cdot 1 + o(h),$$

where the four terms correspond to nothing happening, a birth, a death and a mutation of the original type 1 to type 2. Doing some algebra and letting $h \to 0$

(2.2) $$g_2'(t) = -u_2 g_2(t) - g_2(t)^2 + u_2.$$

If we let $r_1 > r_2$ be the solutions of $x^2 + u_2 x - u_2 = 0$, that is,

(2.3) $$r_i = \frac{-u_2 \pm \sqrt{u_2^2 + 4u_2}}{2},$$

we can write this as

$$g_2'(t) = -(g_2(t) - r_1)(g_2(t) - r_2).$$

Now $g_2(\infty)$ be the probability that a type 2 offspring is eventually generated in the branching process. Letting $t \to \infty$ in (2.2) and noticing that $t \mapsto g_2(t)$ is increasing implies $g_2'(t) \to 0$, we see that

$$0 = -u_2 g_2(\infty) - g_2(\infty)^2 + u_2,$$

so $0 \leq g_2(t) < r_1$ for all $t$ and we have

$$1 = \frac{g_2'(t)}{(r_1 - g_2(t))(g_2(t) - r_2)} = \frac{1}{r_1 - r_2}\left(\frac{g_2'(t)}{g_2(t) - r_2} + \frac{g_2'(t)}{r_1 - g_2(t)}\right).$$

Integrating

$$\ln(g_2(t) - r_2) - \ln(r_1 - g_2(t)) = (r_1 - r_2)t - \ln A,$$

where $A$ is a constant that will be chosen later, so we have

$$\frac{g_2(t) - r_2}{r_1 - g_2(t)} = (1/A)e^{(r_1 - r_2)t}.$$

A little algebra gives

$$g_2(t) = \frac{r_1 + Ar_2 e^{(r_2 - r_1)t}}{1 + Ae^{(r_2 - r_1)t}}.$$

We have $g_2(0) = 0$, so $A = -r_1/r_2$ and

$$g_2(t) = \frac{r_1(1 - e^{(r_2-r_1)t})}{1 - (r_1/r_2)e^{(r_2-r_1)t}}.$$

To prepare for the asymptotics note that (2.3) and the assumption that $u_2 \to 0$ imply that $r_1 - r_2 = \sqrt{u_2^2 + 4u_2} \sim 2\sqrt{u_2}$, $r_1 \sim \sqrt{u_2}$ and $r_1/r_2 \to -1$ so

$$g_2(t) \approx \frac{\sqrt{u_2}(1 - e^{-2\sqrt{u_2}t})}{1 + e^{-2\sqrt{u_2}t}}$$

or to be precise, if $t\sqrt{u_2} \to s$, then

(2.4) $$g_2(t) \sim \sqrt{u_2} \cdot \frac{1 - e^{-2s}}{1 + e^{-2s}}.$$

LEMMA 2.2. *The waiting time for the first type 2 in the branching process with immigration when each type 1 individual experiences mutations at rate $Nu_1$ satisfies*

(2.5) $$Q(\tau_2 \leq t) = 1 - \exp\biggl(-Nu_1 \int_0^t Q_1(\tau_2 \leq s)\,ds\biggr).$$

PROOF. Type 1 mutations are a Poisson process with rate $Nu_1$. A point at time $t - s$ is a success, that is, produces a type 2 before time $t$ with probability $Q_1(\tau_2 \leq s)$. By results for thinning a Poisson process, the number of successes by time $t$ is Poisson with mean $Nu_1 \int_0^t Q_1(\tau_2 \leq s)\,ds$. The result follows from the observation that $Q(\tau_2 \leq t)$ is the probability of at least one success in the Poisson process. □

To find the density function, we recall $g_2(t) = Q_1(\tau_2 \leq t)$ and differentiate to get

$$Nu_1 g_2(t) \exp\biggl(-Nu_1 \int_0^t g_2(s)\,ds\biggr).$$

Changing variables the density function $f_2$ of $\tau_2 \cdot Nu_1\sqrt{u_2}$ is given by

$$f_2(t) = \frac{g_2(t/Nu_1\sqrt{u_2})}{\sqrt{u_2}} \exp\biggl(-Nu_1 \int_0^{t/Nu_1\sqrt{u_2}} g_2(s)\,ds\biggr).$$

Changing variables $r = sNu_1\sqrt{u_2}$ in the integral the above is

$$f_2(t) = \frac{g_2(t/Nu_1\sqrt{u_2})}{\sqrt{u_2}} \exp\biggl(-\int_0^t \frac{g_2(r/Nu_1\sqrt{u_2})}{\sqrt{u_2}}\,dr\biggr).$$

If $Nu_1 \to 0$, then (2.4) implies that the above converges to $\exp(-t)$. If $Nu_1 \to \lambda$, the limit is $h(t)\exp(-\int_0^t h(s)\,ds)$ where

$$h(s) = \frac{1 - e^{-2s/\lambda}}{1 + e^{-2s/\lambda}},$$

which completes the proof of Theorem 1.

**3. A two-type model.** We collect here some results for a simple two-type population model, which we call model $M_0$. We assume that all individuals are either type 0 or type 1, and the population size is always $N$. There are no mutations, and the population evolves according to the Moran model, so each individual dies at rate 1 and then is replaced by a randomly chosen individual in the population. Usually we will assume that the process starts with just one type 1 individual at time zero, but occasionally we will also need to consider starting the process with $j$ type 1 individuals. Denote by $P_j$ and $E_j$ probabilities and expectations when the process is started with $j$ type 1 individuals, and write $P = P_1$ and $E = E_1$. Let $X(t)$ denote the number of type 1 individuals at time $t$.

Let $T_k = \inf\{t: X(t) = k\}$ be the first time at which there are $k$ type 1 individuals, and let $T = \min\{T_0, T_N\}$ be the first time at which all individuals have the same type. Let $L_k$ be the amount of time for which there are $k$ type one individuals, which is the Lebesgue measure of $\{t < T : X(t) = k\}$. Let $R_k$ be the number of times that the number of type 1 individuals jumps to $k$ from $k-1$ or $k+1$. Let $R = 1 + \sum_{k=1}^{N-1} R_k$ be the total number of births and deaths of type 1 individuals. Durrett and Schmidt [8] studied this model and showed that

$$(3.1) \qquad E[R_k | T_0 < T_N] = \frac{2(N-k)^2}{N(N-1)}$$

and

$$(3.2) \qquad E[R_k | T_N < T_0] = \frac{2k(N-k)}{N}.$$

Equation (3.1) is (16) of [8], while (3.2) comes from the beginning of the proof of Lemma 3 in [8].

Because $P(T_N < T_0) = 1/N$, it follows from (3.1) and (3.2) that

$$(3.3) \quad E[R_k] = \frac{(N-1)E[R_k|T_0 < T_N] + E[R_k|T_N < T_0]}{N} = \frac{2(N-k)}{N} \leq 2$$

and, therefore,

$$(3.4) \qquad E[R] = 1 + \sum_{k=1}^{N-1} E[R_k] \leq 2N.$$

If $1 \leq j \leq N-1$, then letting $A$ denote the event that there are at least $j$ type 1 individuals at some time, (3.4) gives

$$(3.5) \qquad E_j[R] = E_j[R\mathbf{1}_A] \leq \frac{E[R\mathbf{1}_A]}{P(A)} = jE[R\mathbf{1}_A] \leq jE[R] \leq 2jN.$$

Turning to the quantities $L_k$, note that when there are $k$ type 1 individuals, births and deaths are each happening at rate $k(N-k)/N$, so the number of type 1 individuals changes again after an exponential time with mean $N/[2k(N-k)]$. Therefore, (3.3) gives

$$(3.6) \qquad E[L_k] = \frac{N}{2k(N-k)} E[R_k] = \frac{1}{k}.$$

Since $P_j(T_k < T_0) = j/k$ for $1 \leq j < N$, we have

$$(3.7) \qquad E_j[L_k] \leq E_1[L_k | T_k < T_0] = \frac{E_1[L_k]}{P_1(T_k < T_0)} = 1,$$

where to emphasize the change in initial condition, we have written $E$ as $E_1$. Since $T = \sum_{k=1}^{N-1} L_k$, it also follows from (3.6) that

$$(3.8) \qquad E[T] = \sum_{k=1}^{N-1} \frac{1}{k} \leq C \log N$$

and it follows from (3.7) that for $j = 1, \ldots, N-1$,

$$(3.9) \qquad E_j[T] \leq N.$$

Finally, we will use branching process theory to obtain the following complement to (3.8).

LEMMA 3.1. *There exists a constant $C$ such that $P(T > t) \leq C/t$ for all $0 \leq t \leq N$.*

PROOF. Consider a continuous-time branching process started with one individual in which each individual dies at rate one and gives birth at rate one. Let $T'$ be the time at which the process becomes extinct. By a theorem of Kolmogorov [16], proved in Section I.9 of [3], and the fact that a Markovian continuous-time branching process can be reduced to a discrete time Galton–Watson process by only examining it at integer times, we see that there is a constant $C'$ such that $P(T' > t) \leq C'/t$ for all $t \geq 0$.

When there are $k$ individuals in the branching process, births and deaths happen at rate $k$. When there are $k$ individuals in the model $M_0$, births and deaths happen at rate $k(N-k)/N$, which is at least $k/2$ as long as $k \leq N/2$. Since the probability that the number of individuals in model $M_0$ ever exceeds $N/2$ is at most $2/N$, we have $P(T > t) \leq 2C'/t + 2/N$ for all $t$, which implies the result. □

**4. The probability of a type $m$ descendant.** We now consider model $M_1$, which evolves in the same way as the process described in the Introduction except that initially there is one type 1 individual and $N-1$ type 0 individuals, and no further type 1 mutations occur. The number of individuals of nonzero type in model $M_1$ therefore evolves exactly like the number of type 1 individuals in model $M_0$, defined at the beginning of the previous section, but in model $M_1$ mutations to types greater than one are possible. The probability, which we denote by $q_m$, that a type $m$ individual is eventually born in model $M_1$ is the same as the probability that a given type one individual in the process described in the Introduction has a type $m$ descendant. Our main goal in this section is to prove the following result.

PROPOSITION 4.1. *Fix an integer $m \geq 2$. Assume conditions* (ii), (iii) *and* (iv) *of Theorem 2 hold. Then $q_m \sim r_{1,m}$.*

We will use Proposition 4.1 to prove Theorem 2. To prove Theorem 3, we will need the following corollary. Here we denote by $q_{j,m}$ the probability that a type $m$ individual eventually appears in a process with initially one type $j$ individual, $N-1$ type 0 individuals, and mutations to type 1 are not allowed.

COROLLARY 4.1. *Fix an integer $m \geq 2$. Assume conditions* (ii) *and* (iii) *of Theorem 2 hold and that $(Nr_{1,m})^2 \to \gamma > 0$. Then $q_{2,m} \sim r_{2,m}$.*

PROOF. We apply the $m-1$ case of Proposition 4.1, with $u_3, \ldots, u_m$ in place of $u_2, \ldots, u_{m-1}$. Since we are assuming (ii) and (iii), we need only to show that $Nr_{2,m} \to \infty$. However, (ii) and (iii) imply

$$\frac{Nr_{2,m}}{Nr_{1,m}} = \frac{u_3^{1/2} u_4^{1/4} \cdots u_m^{1/2^{m-2}}}{u_2^{1/2} u_3^{1/4} \cdots u_{m-1}^{1/2^{m-2}} u_m^{1/2^{m-1}}} > b_2^{1/2} b_3^{1/4} \cdots b_{m-1}^{1/2^{m-2}} u_m^{-1/2^{m-1}} \to \infty.$$

This result and the assumption $(Nr_{1,m})^2 \to \gamma > 0$ imply $Nr_{2,m} \to \infty$. □

We will prove Proposition 4.1 using a branching process approximation. We will approximate model $M_1$ by a continuous-time multi-type branching process in which individuals of type $1 \leq j < m$ die at rate 1, give birth at rate 1 and mutate to individuals of type $j+1$ at rate $u_{j+1}$. Let $p_{j,m}$ be the probability that a type $j$ individual eventually has a descendant of type $m$ in the branching process and let $p_m = p_{1,m}$.

LEMMA 4.1. *If conditions* (ii) *and* (iii) *of Theorem 2 hold, then $p_{j,m} \sim r_{j,m}$.*

PROOF. We proceed by induction starting at $j = m$ and working down to $j = 1$. Clearly, $p_{m,m} = 1$, so the result is valid for $j = m$. Now assume the result is true for $j + 1$. By conditioning on the first event in the branching process, it follows that

$$p_{j,m} = \frac{1}{2 + u_{j+1}}(2p_{j,m} - p_{j,m}^2) + \frac{u_{j+1}}{2 + u_{j+1}}p_{j+1,m}.$$

Multiplying by $2 + u_{j+1}$ and rearranging, we get $p_{j,m}^2 + bp_{j,m} - u_{j+1}p_{j+1,m} = 0$, where $b = u_{j+1}$. The only positive solution is

$$(4.1) \qquad p_{j,m} = \frac{-b + \sqrt{b^2 + 4u_{j+1}p_{j+1,m}}}{2}.$$

Calculus tells that for $h > 0$

$$\sqrt{x+h} - \sqrt{x} = \int_x^{x+h} \frac{1}{2\sqrt{y}}\,dy \leq \frac{h}{2\sqrt{x}},$$

so, we have

$$(4.2) \qquad \begin{aligned} 2\sqrt{u_{j+1}p_{j+1,m}} &\leq \sqrt{4u_{j+1}p_{j+1,m} + b^2} \\ &\leq 2\sqrt{u_{j+1}p_{j+1,m}} + \frac{b^2}{4\sqrt{u_{j+1}p_{j+1,m}}}. \end{aligned}$$

Conditions (ii) and (iii) imply that $u_{j+1} \ll r_{j+1,m}$ and, therefore, that $\sqrt{u_{j+1}r_{j+1,m}} \gg b = u_{j+1}$. Since $p_{j+1,m} \sim r_{j+1,m}$ by the induction hypothesis, it follows from (4.1) and (4.2) that $p_{j,m} \sim \sqrt{u_{j+1}r_{j+1,m}}$. The lemma follows by induction. □

REMARK. One gets the same result for a number of other variants of the model. We leave it to the reader to check that Lemma 4.1 holds when mutation only occurs at birth. To prepare for the proof of Lemma 4.7, we will now show that it holds when type $j$'s give birth to type $j$'s at rate one and to type $j+1$'s at rate $u_{j+1}$. In this case, the first equation is

$$p_{j,m} = \frac{1}{2 + u_{j+1}}(2p_{j,m} - p_{j,m}^2) + \frac{u_{j+1}}{2 + u_{j+1}}(p_{j,m} + p_{j+1,m} - p_{j,m}p_{j+1,m})$$

and rearranges to become $p_{j,m}^2 + u_{j+1}p_{j+1,m}p_{j,m} - u_{j+1}p_{j+1,m} = 0$. Taking $b = u_{j+1}p_{j+1,m}$, the proof goes as before.

We will now prove Proposition 4.1 by induction. We begin with the case $m = 2$, in which the comparison with the branching process is straightforward.

LEMMA 4.2. *Under the assumptions of Proposition 4.1 with $m=2$, we have $q_2 \sim r_{1,2} = u_2^{1/2}$.*

PROOF. If we track the number of type 1 individuals in model $M_1$ before the first type 2 mutation occurs, upward and downward jumps occur at the same rate, which if there are $k$ type 1 individuals is $k(N-k)/N$. For the branching process, when there are $k$ type 1 individuals, upward and downward jumps occur at rate $k$. Therefore, the embedded jump chain (which gives the sequence of states visited by the continuous-time chain) is a simple random walk $S_n$ with $S_0 = 1$ both for model $M_1$ and for the branching process. Therefore, writing $p_2$ as a function of the underlying mutation rate, we claim that for any $L$,

$$(4.3) \qquad p_2(u_2) - 1/N \leq q_2 \leq p_2(u_2 N/(N-L)) + 1/L.$$

The first inequality follows from the fact that unless the number of type 1 individuals in model $M_1$ reaches $N$, which happens with probability $1/N$, model $M_1$ has the same embedded jump chain as the branching process and jumps more slowly. For the second inequality, we note that the probability the Moran model reaches height $L$ is $1/L$. When this does not occur, the Moran model always jumps at rate at least $(N-L)/N$ times the branching process rate. Lemma 4.1 gives $p_2(u_2) \sim u_2^{1/2}$. Condition (iv) gives $Nu_2^{1/2} \to \infty$, so we can choose $L$ such that $L/N \to 0$ and $Lu_2^{1/2} \to \infty$. Under these conditions, (4.3) implies $q_2 \sim u_2^{1/2}$. □

For the rest of this section, we will fix $m$ and assume that the assumptions of Proposition 4.1 hold. We will also assume that Proposition 4.1 has been established for $m-1$, which implies that $q_{2,m} \sim r_{2,m}$. We will reduce the general case to the $m=2$ case in which type 2 mutations occur at rate $u_2 r_{2,m}$. The next two lemmas will allow us to ignore certain type 2 mutations.

LEMMA 4.3. *Let $A_m$ be the event that in model $M_1$ some type 2 mutation that occurs while there is another individual in the population of type 2 or higher has a type $m$ descendant. Then $P(A_m) \ll r_{1,m}$.*

PROOF. Let $\varepsilon > 0$. Let $B$ be the event that the number of individuals in the population of type 1 or higher never exceeds $\varepsilon^{-1} r_{1,m}^{-1}$, so $P(B^c) \leq \varepsilon r_{1,m}$. Let $U = \{t: \text{there is an individual of type 2 or higher alive at time } t\}$. On $B$, type 2 mutations occur at rate at most $\varepsilon^{-1} r_{1,m}^{-1} u_2$ and have a type $m$ descendant with probability $q_{2,m}$. Therefore, letting $|U|$ denote the Lebesgue measure of $U$, we have

$$P(A_m) \leq \varepsilon r_{1,m} + E[|U|\mathbf{1}_B]\varepsilon^{-1} r_{1,m}^{-1} u_2 q_{2,m}.$$

For $k \leq \varepsilon^{-1} r_{1,m}^{-1}$, it follows from (3.6) that the expected amount of time for which there are $k$ individuals of type 1 or higher is $1/k$, and so the expected number of type 2 mutations during this time is at most $(1/k)(ku_2) = u_2$. Therefore, the expected number of type 2 mutations while there are at most $\varepsilon^{-1} r_{1,m}^{-1}$ individuals of type 1 or higher is at most $\varepsilon^{-1} r_{1,m}^{-1} u_2$. By (3.8), the expected amount of time for which these mutations or their offspring are alive in the population is at most $(C \log N) \varepsilon^{-1} r_{1,m}^{-1} u_2$. Therefore, $E[|U|\mathbf{1}_B] \leq (C \log N) \varepsilon^{-1} r_{1,m}^{-1} u_2$. Since $q_{2,m} \sim r_{2,m}$ by the induction hypothesis and $u_2 r_{2,m} = r_{1,m}^2$, it follows that there exists a constant $C$ such that

$$P(A_m) \leq \varepsilon r_{1,m} + C(\log N) \varepsilon^{-2} r_{1,m}^{-2} u_2^2 r_{2,m} = \varepsilon r_{1,m} + C(\log N) \varepsilon^{-2} u_2.$$

Conditions (ii) and (iii) imply that there exist constants $C_1$ and $C_2$ such that

$$\frac{(\log N) u_2}{r_{1,m}} \leq C_1 u_2^{1/2^{m-1}} \log N \leq C_2 u_m^{1/2^{m-1}} \log N \to 0.$$

It follows that

$$\limsup_{N \to \infty} r_{1,m}^{-1} P(A_m) \leq \varepsilon,$$

which implies the lemma. □

LEMMA 4.4. *Let $\varepsilon > 0$. Let $B_m$ be the event that in model $M_1$ some type 2 mutation that occurs while there are fewer than $\varepsilon r_{1,m}^{-1}$ individuals in the population of type 1 or higher has a type $m$ descendant. Then there is a constant $C$, not depending on $\varepsilon$, such that $P(B_m) \leq C \varepsilon r_{1,m}$.*

PROOF. As noted in the proof of Lemma 4.3, the expected number of type 2 mutations while there are $k$ individuals of type 1 or higher is $u_2$. Therefore, the expected number of type 2 mutations while there are fewer than $\varepsilon r_{1,m}^{-1}$ individuals of type 1 or higher is at most $\varepsilon r_{1,m}^{-1} u_2$. By the induction hypothesis, each such mutation produces a type $m$ descendant with probability $q_m \sim r_{2,m}$, so the probability that one of these mutations produces a type 2 descendant is at most $C \varepsilon r_{1,m}^{-1} r_{2,m} u_2$. The desired result now follows from the fact that $u_2 r_{2,m} = r_{1,m}^2$. □

Our strategy is to show that we can reduce the problem to the $m = 2$ case by assuming that each type 2 mutation independently generates a type $m$ descendant with probability $q_{2,m}$. Complicating this picture is the fact that the evolution of the number of type 1 individuals (which produce the type 2 mutations) is not independent of the success of the type 2 mutations because a new individual of type $j \geq 2$ may replace an existing type 1 individual and vice versa. To show that this is not a significant problem, we will construct a coupling of model $M_1$ with another process in which this dependence has been eliminated. We first define model $M_2$ to evolve like model $M_1$ except that initially there are $k$ individuals of type 1 and $N - k$ of type 0, and type 2 mutations are only permitted when there are no individuals of type $j \geq 2$. We then compare model $M_2$ to model $N_2$, in which the type 1 individuals are decoupled from type 2 individuals and their offspring by declaring that (provided a type 0 individual exists):

- if a proposed move exchanges a type 1 and a type $j \geq 2$, we instead exchange a type 0 and a type $j$;
- a mutation that occurs to a type 1 produces a new type 2 individual but replaces a type 0 individual instead of the type 1 that mutated.

To define the coupling precisely, introduce a Poisson process with rate $N$ at which the successive exchanges will occur and let $i_n$ and $j_n$ be independent i.i.d. uniform on $\{1, 2, \ldots, N\}$. In both models, we replace individual $i_n$ with a copy of individual $j_n$. In model $N_2$, if $i_n$ has type 1 and $j_n$ has type 2 or higher, then we choose a type 0 individual at random to become type 1, so that the number of type 1 individuals stays the same. Likewise, if $i_n$ has type 2 or higher and $j_n$ has type 1, then we choose a type 1 individual to become type 0 in model $N_2$. This recipe breaks down when there are no individuals of type 0. However, Lemma 4.5 shows that with high probability the number of individuals of nonzero type is $o(N)$ up to time $\tau_m$. For the mutations, we have for each $1 \leq i \leq N$ a Poisson process with rate $u_2$, which in both models causes a mutation of the $i$th individual, unless either the $i$th individual has type 0 or the $i$th individual has type 1 and there is an individual of type 2 or higher in the population. In model $N_2$, if a type 1 individual mutates to type 2, a type 0 individual is chosen at random to become type 1, to keep the number of type 1 individuals constant.

Let $X_1(t)$ and $Y_1(t)$ be the number of type 1 individuals at time $t$ in models $M_2$ and $N_2$, respectively. Let $Z(t) = X_1(t) - Y_1(t)$. Let $\hat{X}_2(t)$ and $\hat{Y}_2(t)$ be the number of individuals in models $M_2$ and $N_2$, respectively, of type greater than or equal to 2. Note that by renumbering the individuals as the process evolves if necessary, we can ensure that for all $t \geq 0$, at time $t$ there are $\min\{X_1(t), Y_1(t)\}$ integers $j$ such that the $j$th individual has type 1 in both model $M_2$ and model $N_2$. Note also that with the above coupling, if a type 2 mutation occurs at the same time in both models, descendants of this mutation will always have the same type in both models. This means that if the mutation has a type $m$ descendant in one model, then it will have a type $m$ descendant in the other. Finally, as long as the number of individuals of nonzero type stays below $N/2$, we can also ensure that there is no $j$ such that the $j$th individual has type 1 in one of the two models and type 2 or higher in the other. The lemma below, combined with condition (iv), ensures that in both models, the number of individuals of nonzero type stays much smaller than $N$.

LEMMA 4.5. *Fix $t > 0$. Suppose $X_1(0) = Y_1(0) = [\varepsilon r_{1,m}^{-1}]$ and $\hat{X}_2(0) = \hat{Y}_2(0) = 0$. Assume $f$ is a function of $N$ such that $f(N) r_{1,m} \to \infty$ as $N \to \infty$. Then using $\to_p$ to denote convergence in probability, we have*

$$\max_{0 \leq s \leq tr_{1,m}^{-1}} \frac{X_1(s) + \hat{X}_2(s)}{f(N)} \to_p 0 \quad \text{and} \quad \max_{0 \leq s \leq tr_{1,m}^{-1}} \frac{Y_1(s) + \hat{Y}_2(s)}{f(N)} \to_p 0.$$

PROOF. In model $M_2$, individuals of type 1 or higher give birth and die at the same rate, so $(X_1(s) + \hat{X}_2(s), s \geq 0)$ is a martingale and

$$E[X_1(tr_{1,m}^{-1}) + \hat{X}_2(tr_{1,m}^{-1})] = X_1(0) + \hat{X}_2(0) = [\varepsilon r_{1,m}^{-1}].$$

By Doob's maximal inequality, if $\delta > 0$, then

$$P\left(\max_{0 \leq s \leq tr_{1,m}^{-1}} \frac{X_1(s) + \hat{X}_2(s)}{f(N)} > \delta\right) \leq \frac{E[X_1(tr_{1,m}^{-1}) + \hat{X}_2(tr_{1,m}^{-1})]}{\delta f(N)}$$

$$\leq \frac{\varepsilon r_{1,m}^{-1}}{\delta f(N)} \to 0$$

as $N \to \infty$, which implies the first statement of the lemma.

In model $N_2$, mutations of type 1 individuals cause new type 2 individuals to replace type 0 individuals. Births and deaths occur at the same rate, so the process $(Y_1(s), s \geq 0)$ is a martingale, while $(Y_1(s) + \hat{Y}_2(s), s \geq 0)$ is a submartingale. Now $E[Y_1(s)] = [\varepsilon r_{1,m}^{-1}]$ for all $s$, so the expected number of type 2 individuals that appear before time $tr_{1,m}^{-1}$ because of mutation is at most $\varepsilon r_{1,m} \cdot tr_{1,m}^{-1} \cdot u_2 = \varepsilon u_2 r_{1,m}^{-2} t$. It follows that

$$E[Y_1(tr_{1,m}^{-1}) + \hat{Y}_2(tr_{1,m}^{-1})] \leq \varepsilon r_{1,m}^{-1} + \varepsilon u_2 r_{1,m}^{-2} t.$$

Now

$$(4.4) \quad u_2 r_{1,m}^{-1} = \frac{u_2}{u_2^{1/2} u_2^{1/4} \cdots u_m^{1/2^{m-1}}} = \frac{u_2^{1-1/2^{m-1}}}{u_2^{1/2} u_2^{1/4} \cdots u_m^{1/2^{m-1}}} \cdot u_2^{1/2^{m-1}} \to 0,$$

because condition (ii) implies that the first factor is bounded by a constant, so Doob's maximal inequality this time gives

$$P\left(\max_{0\leq s\leq tr_{1,m}^{-1}} \frac{Y_1(s) + \hat{Y}_2(s)}{f(N)} > \delta\right) \leq \frac{\varepsilon r_{1,m}^{-1} + \varepsilon u_2 r_{1,m}^{-2} t}{\delta f(N)} \to 0,$$

which implies the second half of the lemma. □

We now work on bounding the process $(Z(t), t \geq 0)$. There are three types of events that cause this process to jump. First, whenever a type 1 individual in model $M_2$ mutates to type 2, there is no corresponding change in model $N_2$, because any new type 2 individual in model $N_2$ resulting from mutation replaces a type 0. These changes cause the $Z$ process to decrease by one. Letting $\mu(t)$ be the rate at which they are occurring at time $t$, we have

$$0 \leq \mu(t) \leq u_2 X_1(t),$$

where the second inequality could be strict because mutations are suppressed if there is already a type 2 individual in the population.

Second, one of the "extra" $|Z(t)|$ type 1 individuals in one process or the other could experience a birth or a death. This could cause the $Z$ process to increase or decrease by one. If $X_1(t) > Y_1(t)$, then at time $t$, both increases and decreases in the $Z$ process occur because of such changes at rate $|Z(t)|(N - |Z(t)|)/N$, because the $Z$ process changes unless the other individual involved in the exchange was also one of the $|Z(t)|$ individuals that are type 1 in model $M_2$ but not model $N_2$. If $Y_1(t) > X_1(t)$, then increases and decreases in the $Z$ process occur at rate $|Z(t)|(N - |Z(t)| - \hat{Y}_2(t))/N$ because exchanges between a type 1 individual and an individual of type 2 or higher are prohibited in model $N_2$.

Finally, there are transitions in which one of the $\min\{X_1(t), Y_1(t)\}$ individuals that are type 1 in both processes experiences a birth or death, but the other individual involved in the exchange is one of the $\hat{Y}_2(t)$ individuals that has type 2 in model $N_2$, so the type 1 population does not change in model $N_2$. Such changes occur at rate $\hat{Y}_2(t) \min\{X_1(t), Y_1(t)\}/N$.

Thus, if we let

$$\lambda(t) = \frac{|Z(t)|(N - |Z(t)| - \hat{Y}_2(t)\mathbf{1}_{\{Y_1(t) > Z_1(t)\}})}{N} + \frac{\hat{Y}_2(t)\min\{X_1(t), Y_1(t)\}}{N},$$

then at time $t$ the $Z$ process is increasing by 1 at rate $\lambda(t)$ and decreasing by 1 at rate $\lambda(t) + \mu(t)$. The next result uses these facts to control the difference between $X_1(t)$ and $Y_1(t)$.

LEMMA 4.6. *Fix $t > 0$. Let $Z_N(s) = r_{1,m} Z(sr_{1,m}^{-1})$ for all $s \geq 0$. If $X_1(0) = Y_1(0) = \varepsilon r_{1,m}^{-1}$ and $\hat{X}_2(0) = \hat{Y}_2(0) = 0$, then*

$$\max_{0\leq s\leq t} Z_N(s) \to_p 0.$$

PROOF. We will use Theorem 4.1 from Chapter 7 in [10] to show that $Z_N$ converges to a diffusion with $b(x) = 0$, $a(x) = 2|x|$, and initial point 0, so the limit is identically zero. The first step is to observe that the Yamada–Watanabe theorem; see, for example, (3.3) on page 193 of [6], gives pathwise uniqueness for the limiting SDE, which in turn implies that the martingale problem is well posed. To verify the other assumptions of the theorem, define

$$B_N(t) = -\int_0^t \mu(sr_{1,m}^{-1})\, ds$$

and

$$A_N(t) = \int_0^t r_{1,m}(2\lambda(sr_{1,m}^{-1}) + \mu(sr_{1,m}^{-1}))\, ds.$$

In view of the transition rates for the process $(Z(t), t \geq 0)$, we see that at time $s$ the process $Z_N(s)$ experiences positive jumps by the amount $r_{1,m}$ at rate $\lambda(sr_{1,m}^{-1})r_{1,m}^{-1}$ and negative jumps by the same amount at rate $(\lambda(sr_{1,m}^{-1}) + \mu(sr_{1,m}))r_{1,m}^{-1}$. Therefore, letting $M_N(t) = Z_N(t) - B_N(t)$, the processes $(M_N(t), t \geq 0)$ and $(M_N^2(t) - A_N(t), t \geq 0)$ are martingales. To obtain the result of the lemma from Theorem 4.1 in Chapter 7 of [10], it remains to show that for any fixed $T > 0$, we have

(4.5) $$\sup_{0 \leq t \leq T} |B_N(t)| \to_p 0$$

and

(4.6) $$\sup_{0 \leq t \leq T} \left| A_N(t) - \int_0^t 2|Z_N(s)|\, ds \right| \to_p 0.$$

To prove (4.5), note that

$$\sup_{0 \leq t \leq T} |B_N(t)| \leq T \sup_{0 \leq t \leq T} \mu(tr_{1,m}^{-1}) \leq Tu_2 \max_{0 \leq t \leq Tr_{1,m}^{-1}} X_1(t).$$

Since $r_{1,m}/(Tu_2) \to \infty$ by (4.4), (4.5) now follows from Lemma 4.5 with $f(N) = 1/(Tu_2)$. For (4.6), note that

$$A_N(t) - \int_0^t 2|Z_N(s)|\, ds$$
$$= r_{1,m} \int_0^t \Bigg( -\frac{2|Z(sr_{1,m}^{-1})|^2}{N} - \frac{2|Z(sr_{1,m}^{-1})|\hat{Y}_2(sr_{1,m}^{-1})\mathbf{1}_{\{Y_1(sr_{1,m}^{-1}) > Z_1(sr_{1,m}^{-1})\}}}{N}$$
$$+ \frac{2\hat{Y}_2(sr_{1,m}^{-1})\min\{X_1(sr_{1,m}^{-1}), Y_1(sr_{1,m}^{-1})\}}{N} + \mu(sr_{1,m}^{-1}) \Bigg) ds.$$

It suffices to control the absolute values of the four terms over all $t \leq T$. Note that $Z(sr_{1,m}^{-1}) \leq \max\{X_1(sr_{1,m}^{-1}), Y_1(sr_{1,m}^{-1})\}$. Therefore, by Lemma 4.5 with $f(N) = \sqrt{N/r_{1,m}}$, the three quantities $\max_{0 \le s \le Tr_{1,m}^{-1}} r_{1,m}^{1/2} N^{-1/2} |Z(s)|$, $\max_{0 \le s \le Tr_{1,m}^{-1}} r_{1,m}^{1/2} N^{-1/2} \hat{Y}_2(s)$ and $\max_{0 \le s \le Tr_{1,m}^{-1}} r_{1,m}^{1/2} N^{-1/2} X_1(s)$ all converge in probability to zero as $N \to \infty$. This is enough to establish the convergence of the first three terms. The result for the third term follows from (4.5) and the fact that $r_{1,m} \to 0$. $\square$

In the model $N_2$, types $j \ge 2$ have the same relationship to type 1 individuals as in the branching process. That is, type 1's give birth to type 2's, but the fate of a type 2 family does not affect the number of type 1 individuals because a type 1 individual cannot be exchanged with an individual of type 2 or higher. Lemma 4.3 has shown that we can ignore type 2 births that occur when another type 2 is present, so successive type 2 births give independent chances of producing a type $m$ individual. We are now close to our goal announced in the Introduction of reducing the $m$-type problem to the 2-type problem with $\bar{u}_2 = u_2 q_{2,m}$, that is, to the simplified model in which at each type 2 mutation, we flip a coin with probability $q_{2,m}$ of heads to see if it will generate a type $m$ individual.

Let model $N_2'$ be model $N_2$ modified so that if a type 2 mutation occurs when $\hat{Y}_2(t) > 0$, instead of suppressing this event entirely, we flip a coin with probability $q_{2,m}$ of heads. We then add a type $m$ individual to the population if the coin is heads and otherwise make no change. Lemma 4.3 implies that the difference between the probability of getting a type $m$ individual in model $N_2$ and the probability of getting a type $m$ individual in model $N_2'$ tends to zero as $N \to \infty$. However, it is easier to prove the next result using model $N_2'$ because in model $N_2'$, each type 1 individual is giving rise to individuals that will produce a type $m$ descendant at rate $u_2 q_{2,m}$, regardless of whether there are other individuals in the population of type 2 or higher.

LEMMA 4.7. *Let $\varepsilon > 0$. Consider model $N_2'$ starting from $[\varepsilon r_{1,m}^{-1}]$ type 1 individuals at time zero. Let $h_{N,m,\varepsilon}^1$ be the probability that a type $m$ individual is born at some time. Then*
$$\lim_{N \to \infty} h_{N,m,\varepsilon}^1 = 1 - e^{-\varepsilon}.$$

PROOF. Consider a modified branching process in which type $j$ individuals give birth at rate one, die at rate one, and give birth to type $j+1$ individuals at rate $u_{j+1}$. Let $h_{N,m,\varepsilon}^0$ be the probability that if the branching process starts with $[\varepsilon r_{1,m}]$ individuals, a type $m$ individual is born at some time. Since different families are independent, Lemma 4.1 implies
$$h_{N,m,\varepsilon}^0 = 1 - (1-p_m)^{[\varepsilon r_{1,m}^{-1}]} \to 1 - e^{-\varepsilon},$$
where $p_m$ is the probability that a type 1 individual has a type $m$ descendant.

We now compare this process to model $N_2'$. The number of type 1 individuals in model $N_2'$ jumps more slowly than the number of type 1 individuals in the branching process, but in both processes type 1 individuals give birth to type 2 individuals at rate $u_2$, and then type 2 individuals and their descendants evolve independently of the type 1's. Therefore, if the probability $p_{2,m}$ that a type 2 individual in the branching process produces a type $m$ descendant were equal to $q_{2,m}$, then it would follow that $h^1_{N,m,\varepsilon} \geq h^0_{N,m,\varepsilon}$. Instead, we only have $p_{2,m} \sim q_{2,m}$ because $p_{2,m} \sim r_{2,m}$ by the remark after Lemma 4.1 and $q_{2,m} \sim r_{2,m}$ by the induction hypothesis. It follows that

$$h^1_{N,m,\varepsilon} \geq h^0_{N,m,\varepsilon}(1 - o(1)) \to 1 - e^{-\varepsilon}.$$

To get a bound in the opposite direction, observe that we can pick $K \to \infty$ so that $L = Kr^{-1}_{1,m} = o(N)$, and with probability tending to one as $N \to \infty$, the number of type 1's does not reach $L$. Therefore, writing $h^1_{N,m,\varepsilon}$ and $h^0_{N,m,\varepsilon}$ as functions of the rate at which type 1 individuals give birth to type 2 individuals, we have

$$h^1_{N,m,\varepsilon}(u_2) \leq h^0_{N,m,\varepsilon}(u_2 N/(N-L))(1 + o(1)) + o(1) \to 1 - e^{-\varepsilon},$$

which completes the proof. □

LEMMA 4.8. *Let $\varepsilon > 0$. Consider model $M_2$ starting from $[\varepsilon r^{-1}_{1,m}]$ type 1 individuals at time zero. Let $h_{N,m,\varepsilon}$ be the probability that a type $m$ individual is born at some time. Then*

$$\lim_{N \to \infty} |h_{N,m,\varepsilon} - h^1_{N,m,\varepsilon}| = 0.$$

PROOF. Recall the coupling between model $M_2$ and model $N_2$ described earlier in this section. With this coupling, if a type 2 mutation occurs at the same time in both processes, then it produces a type $m$ descendant in one process if and only if it produces a type $m$ descendant in the other. Consequently, it suffices to bound the probability that some type 2 mutation that appears in one process but not the other produces a type $m$ descendant.

Lemma 4.3 bounds this probability for mutations that occur in one model but get suppressed in the other because there are no individuals of type 2 or higher. It remains to consider the mutations experienced by the $|Z(t)|$ individuals that are type 1 in one process but not the other. Pick $s$ large enough so that the probability $N_2$ or $M_2$ does not die out by time $sr^{-1}_{1,m}$ is $< \delta$. Pick $\eta$ so that $\eta s < \delta^2$. By Lemma 4.6, if $N$ is large, we have $\max_{t \leq s} |Z_N(t)| < \eta$ with probability $> 1 - \delta$. The expected number of births that occur in one process but not in the other before time $sr^{-1}_{1,m}$ when $\max_{t \leq s} |Z_N(t)| < \eta$ is bounded by

$$2\eta r^{-1}_{1,m} \cdot sr^{-1}_{1,m} u_2 \leq 2\delta^2 r^{-2}_{1,m} u_2.$$

Using Chebyshev's inequality, it follows that with probability $> 1 - 4\delta$ the number of type 2 mutant births that occur in one process but not the other is bounded by $\delta r_{1,m}^{-2} u_2 = \delta r_{2,m}^{-1}$. When this occurs, the success probabilities differ by at most $\delta$ because each mutation has probability $q_{2,m} \sim r_{2,m}$ of producing a type $m$ descendant. Since $\delta > 0$ is arbitrary, the desired results follow. $\square$

PROOF OF PROPOSITION 4.1. The probability that the number of individuals of type greater than zero reaches $[\varepsilon r_{1,m}^{-1}]$ is $1/[\varepsilon r_{1,m}^{-1}]$. If, at the time $T$ when the number of individuals of type greater than zero reaches $[\varepsilon r_{1,m}^{-1}]$, we change the type of all individuals whose type is nonzero to type 1, and if we disregard type 2 mutations that occur when there is another individual of type $j \geq 2$, then the probability of getting a type $m$ individual after this time becomes $h_{N,m,\varepsilon}$. Since these changes of the types can only reduce the probability of getting a type $m$ individual, we have

$$(4.7) \qquad q_m \geq \frac{1}{[\varepsilon r_{1,m}^{-1}]} h_{N,m,\varepsilon}.$$

Also, for a type $m$ individual to appear, either the type $m$ individual must be descended from a type 1 individual that is alive at time $T$, or else the type $m$ individual must be descended from a type 2 individual that existed before time $T$, so using Lemmas 4.3 and 4.4, it follows that

$$(4.8) \qquad q_m \leq \frac{1}{[\varepsilon r_{1,m}^{-1}]} h_{N,m,\varepsilon} + C\varepsilon r_{1,m}.$$

The result follows by letting $\varepsilon \to 0$. $\square$

**5. Proof of Theorem 2.** In this section, we complete the proof of Theorem 2. The argument is based on the following result on Poisson approximation, which is part of Theorem 1 of [2].

LEMMA 5.1. *Suppose $(A_i)_{i \in \mathcal{I}}$ is a collection of events, where $\mathcal{I}$ is any index set. Let $W = \sum_{i \in \mathcal{I}} \mathbf{1}_{A_i}$ be the number of events that occur, and let $\lambda = E[W] = \sum_{i \in \mathcal{I}} P(A_i)$. Suppose for each $i \in \mathcal{I}$, we have $i \in B_i \subset \mathcal{I}$. Let $\mathcal{F}_i = \sigma((A_j)_{j \in \mathcal{I} \setminus B_i})$. Define*

$$b_1 = \sum_{i \in \mathcal{I}} \sum_{j \in B_i} P(A_i) P(A_j),$$

$$b_2 = \sum_{i \in \mathcal{I}} \sum_{i \neq j \in B_i} P(A_i \cap A_j),$$

$$b_3 = \sum_{i \in \mathcal{I}} E[|P(A_i | \mathcal{F}_i) - P(A_i)|].$$

*Then $|P(W = 0) - e^{-\lambda}| \leq b_1 + b_2 + b_3$.*

We will use the following lemma to get the second moment estimate needed to bound $b_2$. When we apply this result, the individuals born at times $t_1$ and $t_2$ will both be type 1. We call the second one type 2 to be able to easily distinguish the descendants of the two individuals.

LEMMA 5.2. *Fix times $t_1 < t_2$. Consider a population of size $N$ which evolves according to the Moran model in which all individuals initially have type 0. There are no mutations, except that one individual becomes type 1 at time $t_1$, and one type 0 individual (if there is one) becomes type 2 at time $t_2$. Fix a positive integer $L \leq N/2$. For $i = 1, 2$, let $Y_i(t)$ be the number of type $i$ individuals at time $t$ and let $B_i$ be the event that $L \leq \max_{t \geq 0} Y_i(t) \leq N/2$. Then*

$$P(B_1 \cap B_2) \leq 2/L^2.$$

PROOF. Because $(Y_1(t), t \geq t_1)$ is a martingale, it is clear that $P(B_1) \leq 1/L$. Let $s_1 < s_2 < \cdots < s_J$ be the ordered times, after time $t_2$, at which the $Y_1$ process jumps. For $t \geq t_2$, let $Z(t) = Y_2(t)A(t)$, where

$$A(t) = \frac{N - Y_1(t_2)}{N - Y_1(t)} = \prod_{i\,:\,s_i \leq t} \frac{N - Y_1(s_i-)}{N - Y_1(s_i)}.$$

We claim that conditional on $(Y_1(t), t \geq t_1)$, the process $(Z(t), t \geq t_2)$ is a martingale.

To see this, note that between the times $s_i$, births and deaths of type 2 individuals occur at the same rate, even conditional on $(Y_1(t), t \geq t_1)$, so $Z(t)$ experiences both positive and negative jumps of size $(N - Y_1(t_2))/(N - Y_1(t))$ at the same rate. At the time $s_i$, if $Y_1(s_i) = Y_1(s_i-) + 1$, then one of the $N - Y_1(s_i-)$ individuals of type other than 1 dies at time $s_i$, so we have $Y_2(s_i) = Y_2(s_i-) - 1$ with probability $\alpha_i = Y_2(s_i-)/(N - Y_1(s_i-))$ and $Y_2(s_i) = Y_2(s_i-)$ with probability $1 - \alpha_i$. Note that

$$(1 - \alpha_i)Y_2(s_i-) + \alpha_i(Y_2(s_i-) - 1) = Y_2(s_i-) - \alpha_i$$
$$= Y_2(s_i-)\left(1 - \frac{1}{N - Y_1(s_i-)}\right)$$
$$= Y_2(s_i-)\frac{N - Y_1(s_i)}{N - Y_1(s_i-)}.$$

Likewise, if $Y_1(s_i) = Y_1(s_i-) - 1$, then one of the $N - Y_1(s_i-)$ individuals of type other than 1 gives birth at time $s_i$, so $Y_2(s_i) = Y_2(s_i-) + 1$ with probability $\alpha_i = Y_2(s_i-)/(N - Y_1(s_i-))$ and $Y_2(s_i) = Y_2(s_i-)$ with probability $1 - \alpha_i$, and we have

$$(1 - \alpha_i)Y_2(s_i-) + \alpha_i(Y_2(s_i-) + 1) = Y_2(s_i-) + \alpha_i$$

$$= Y_2(s_i-)\left(1 + \frac{1}{N - Y_1(s_i-)}\right)$$

$$= Y_2(s_i-)\frac{N - Y_1(s_i)}{N - Y_1(s_i-)}.$$

The martingale property follows because $A(s_i) = A(s_i-)(N-Y_1(s_i-))/(N-Y_1(s_i))$, compensating for the expected change in the $Y_2$ process.

Since $(Z(t), t \geq t_2)$ is a martingale conditional on $(Y_1(t), t \geq t_1)$ and $Z(t_2) = 1$, we have $P(Z(t) \geq L/2$ for some $t|B_1) \leq 2/L$. On the event $B_1$, we have $A(t) \leq 2$ for all $t \geq t_2$, so

$$P(B_2|B_1) \leq P(Y_2(t) \geq L \text{ for some } t|B_1)$$
$$\leq P(Z(t_2) \geq L/2 \text{ for some } t|B_1) \leq 2/L.$$

Since $P(B_1) \leq 1/L$, the result follows. $\square$

We now introduce a set-up that will allow us to apply Lemma 5.1. Let $\varepsilon > 0$, and let $K$ be a large positive number that will be chosen later. Let $\bar{q}_m$ be the probability that in model $M_1$:

- there is eventually a type $m$ individual in the population,
- the maximum number of individuals of nonzero type over all times is between $\varepsilon/r_{1,m}$ and $N/2$, and
- the family lives for time $\leq K/r_{1,m}$; that is, there are no individuals of nonzero type remaining at time $K/r_{1,m}$.

We will call the second and third points the *side conditions*. Divide the interval $[0, t/(Nr_{0,m})]$ into $M$ subintervals of equal length, where $Mr_{1,m} \to \infty$ as $N \to \infty$. Label the intervals $I_1, \ldots, I_M$, and let $D_i$ be the event that there is a type 1 mutation in the interval $I_i$.

For bookkeeping purposes, we will also introduce type 1b mutations, which individuals of type greater than zero experience at rate $u_1$ but which do not affect the type of the individual. Mutations to type zero individuals will be called type 1a mutations, and the phrase "type 1 mutation" will refer both to type 1a and type 1b mutations for the rest of this section. This will ensure that type 1 mutations are always occurring at rate exactly $Nu_1$. To determine whether or not the first type 1b mutation in interval $i$ is "successful," we let $\xi_1, \ldots, \xi_M$ be i.i.d. random variables, independent of our process, that equal 1 with probability $\bar{q}_m$.

Let $A_i$ be the event that there is a type 1 mutation in the interval $I_i$ and one of the following occurs:

- The first type 1 mutation in $I_i$ has type 1a. The individual that gets this mutation has a type $m$ descendant and the side conditions hold. That is, the maximum number of descendants of the mutation over all times is between $\varepsilon/r_{1,m}$ and $N/2$, and there are no descendants of the mutation remaining in the population at the time $K/r_{1,m}$ after the mutation occurred.
- The first type 1 mutation in $I_i$ has type 1b, and $\xi_i = 1$.

As in Lemma 5.1, let $W = \sum_{i=1}^{M} \mathbf{1}_{A_i}$ be the number of events that occur, and let $\lambda = E[W]$.

LEMMA 5.3. $\limsup_{N \to \infty} |P(W = 0) - e^{-\lambda}| = 0$.

PROOF. Let $\beta_i$ consist of all subintervals whose distance to $I_i$ is at most $K/r_{1,m}$. Define $b_1$, $b_2$ and $b_3$ as in Lemma 5.1. We first claim that $b_3 = 0$. Suppose $I_i = [a, b]$. Note that the event $A_i$ does not depend on the state of the population at time $a$. Also, because of the side condition that a mutation is not considered successful if it has descendants surviving for a time longer than $K/r_{1,m}$, the event $A_i$ is determined by time $b + K/r_{1,m}$ and is therefore independent of the events $A_j$ for $j > i$ and $j \notin \beta_i$. Likewise, all of the events $A_j$ for $j < i$ and $j \notin \beta_i$ are determined by the behavior of the process before time $a$, so these events are also independent of $A_i$. It follows that $A_i$ is independent of $(A_j)_{j \notin \beta_i}$, and thus that $b_3 = 0$.

The length $|I_i|$ of the interval $I_i$ is $t/(MNr_{0,m})$, so since type 1 mutations occur at rate $Nu_1$, we have $P(D_i) \leq Nu_1|I_i| = t/(Mr_{1,m})$. Since $P(A_i|D_i) = \bar{q}_m$, Proposition 4.1 gives

$$P(A_i) = \bar{q}_m P(D_i) \leq t\bar{q}_m/(Mr_{1,m}) \sim t/M.$$

There are at most $2(K/(r_{1,m}|I_i|) + 1)$ intervals in $\beta_i$, so for large $M$

$$b_1 \leq M \cdot 2\left(\frac{K}{r_{1,m}|I_i|} + 1\right) \cdot \left(\frac{t}{M}\right)^2$$
$$= 2M \cdot \frac{KMNr_{0,m}}{r_{1,m}t}\left(\frac{t}{M}\right)^2 + \frac{2t^2}{M} = 2KNu_1t + \frac{2t^2}{M}.$$

Since $Nu_1 \to 0$ by (i) and $M \to \infty$, $b_1 \to 0$.

To bound $b_2$, note that $P(D_i \cap D_j) \leq [t/(Mr_{1,m})]^2$ because mutations in disjoint intervals occur independently. We now apply Lemma 5.2 with $L = \varepsilon/r_{1,m}$, $t_1$ being the time of the first mutation in $I_i$, and $t_2$ being the time of the first mutation in $I_j$. For the event $A_i$ to occur, it is necessary for the event $B_i$ considered in Lemma 5.2 to occur. Note that mutations do not effect the result of Lemma 5.2 because the side conditions involve all descendants of the original mutation, regardless of type.

$$P(A_i \cap A_j | D_i \cap D_j) \leq 2r_{1,m}^2/\varepsilon^2$$

and thus $P(A_i \cap A_j) \leq 2t^2/(M\varepsilon)^2$. Since there are at most $2(K/(r_{1,m}|I_i|)+1)$ intervals in $\beta_i$, we have

$$b_2 \leq M \cdot 2 \left( \frac{K}{r_{1,m}|I_i|} + 1 \right) \frac{2t^2}{(M\varepsilon)^2}$$

$$= 4M \cdot \frac{KMNr_{0,m}}{r_{1,m}t} \left( \frac{t}{M\varepsilon} \right)^2 + \frac{4t^2}{M\varepsilon^2} = 4\varepsilon^{-2} KNu_1 t + \frac{4t^2}{M\varepsilon^2}.$$

This shows $b_2 \to 0$, and completes the proof. $\square$

LEMMA 5.4. *Let $\sigma_m$ be the first time at which there is a type 1 individual in the population that will have a type m descendant. Then*

(5.1) $$\lim_{N \to \infty} P(\sigma_m > t/(Nr_{0,m})) = \exp(-t).$$

PROOF. To obtain (5.1) from Lemma 5.3, it suffices to show that there is a constant $C$ such that for sufficiently large $N$, we have $|t - \lambda| \leq C\varepsilon$ and $|P(W = 0) - P(\sigma_m > t/(Nr_{0,m}))| \leq C\varepsilon$. The result will then follow by letting $\varepsilon \to 0$. Clearly, $\bar{q}_m \leq q_m$, and $q_m - \bar{q}_m$ is at most the probability that in model $M_1$, a type $m$ individual appears even though either (a) the total number of individuals of nonzero type never exceeds $\varepsilon r_{1,m}$, (b) the total number of individuals of nonzero type exceeds $N/2$, or (c) the family does not die out before $K/r_{1,m}$. Because $Nr_{1,m} \to \infty$, we have $K/r_{1,m} < N$ for sufficiently large $N$, so we can apply Lemma 3.1 to show that the probability that a given mutation survives for as long as $K/r_{1,m}$ is at most $Cr_{1,m}/K$. Using Lemma 4.4, we get

$$q_m - \bar{q}_m \leq C\varepsilon r_{1,m} + 2/N + Cr_{1,m}/K.$$

Since $Nr_{1,m} \to \infty$ by (iv), we have $2/N \ll r_{1,m}$, so if $K$ is large, we get

(5.2) $$q_m - C\varepsilon r_{1,m} \leq \bar{q}_m \leq q_m.$$

Note that

$$\lambda = \sum_{i \in \mathcal{I}} P(A_i) = \sum_{i \in \mathcal{I}} P(D_i) \bar{q}_m = MP(D_1)\bar{q}_m$$

$$= M\bar{q}_m(1 - e^{-Nu_1|I_1|}) \sim M\bar{q}_m Nu_1|I_1| = t\bar{q}_m/r_{1,m}.$$

Because $q_m \sim r_{1,m}$ by Proposition 4.1, this result combined with (5.2) implies $|t - \lambda| \leq C\varepsilon$ for sufficiently large $N$.

It remains to bound $|P(W = 0) - P(\sigma_m > t/(Nr_{0,m}))|$. We can have $W > 0$ with $\sigma_m > t/(Nr_{0,m})$ only if for some $i$, there is a type 1b mutation in $I_i$ and $\xi_i = 1$. Let $X(t)$ be the number of individuals of nonzero type. As long as $X(t)$ stays below $\varepsilon N$, type 1b mutations occur at rate at most $N\varepsilon u_1$, so the probability that this occurs is at most

$$(\varepsilon N u_1)(t/Nr_{0,m})\bar{q}_m \leq C\varepsilon,$$

using Proposition 4.1. Since individuals give birth and die at the same rate, $(X(t), t \geq 0)$ is a submartingale. Also, $E[X(t/(Nr_{0,m}))]$ is the expected number of type 1a mutations before time $t/(Nr_{0,m})$, which is at most $t/r_{1,m}$. Therefore, by Doob's maximal inequality,

$$P(X(s) \geq \varepsilon N \text{ for some } s \leq t/(Nr_{0,m})) \leq t/(\varepsilon Nr_{1,m}),$$

which goes to zero as $N \to \infty$ by condition (iv).

We can have $W = 0$ with $\sigma_m \leq t/(Nr_{0,m})$ in one of two ways. One possibility is that there could be a successful type 1 mutation in one of the $M$ subintervals that is not the first type 1 mutation in that interval. The expected number of type 1 mutations in the $i$th interval that are not the first in their interval is at most $(t/Mr_{1,m})^2$. Therefore, the probability that some successful type 1 mutation is not the first type 1 mutation in its interval is at most $M(t/Mr_{1,m})^2 q_m \leq C/(Mr_{1,m})$. Since $Mr_{1,m} \to \infty$, this probability tends to zero as $N \to \infty$. The other possibility is that there could be a successful type 1 mutation that does not satisfy the extra conditions we imposed. The probability that this occurs is at most

$$(Nu_1)(t/Nr_{0,m})(q_m - \bar{q}_m) \leq Ct\varepsilon$$

by (5.2). This observation completes the proof of the lemma. $\square$

The following result in combination with Lemma 5.4 implies Theorem 2.

LEMMA 5.5. *We have*

(5.3) $$Nr_{0,m}(\tau_m - \sigma_m) \to 0 \quad \text{in probability.}$$

PROOF. Let $\varepsilon > 0$ and $\delta > 0$. By Lemma 5.4, we can choose $s$ large enough that for sufficiently large $N$,

$$P(\sigma_m > s/(Nr_{0,m})) < \delta/3.$$

By Lemma 3.1, the probability that a type 1a mutation takes longer than time $\varepsilon/(Nr_{0,m})$ to die out or fixate is at most $C \max\{1/N, Nr_{0,m}/\varepsilon\}$. Because the expected number of type 1a mutations before time $s/Nr_{0,m}$ is at most $(Nu_1)(s/Nr_{0,m}) = u_1 s/r_{0,m}$, it follows from Markov's inequality that the probability that some type 1a mutation that appears before time $s/(Nr_{0,m})$ takes longer than time $\varepsilon/(Nr_{0,m})$ to die out or fixate is at most $Cs \max\{u_1/(Nr_{0,m}), Nu_1/\varepsilon\}$. As $N \to \infty$, the first of these terms goes to zero by (iv) while the second goes to zero by (i), so this probability is less than $\delta/3$ for sufficiently large $N$. Finally, the probability that one of the type 1a mutations before time $s/(Nr_{0,m})$ fixates is at most

$$\frac{s}{Nr_{0,m}} \cdot Nu_1 \cdot \frac{1}{N},$$

since mutations occur at rate $Nu_1$ and fix with probability $1/N$. This is less than $\delta/3$ for large $N$ by (iv). Hence, $P(Nr_{0,m}(\tau_m - \sigma_m) > \varepsilon) < \delta$ for sufficiently large $N$. □

**6. The key to the proof of Theorem 3.** Throughout this section and the next, we assume all of the hypotheses of Theorem 3 are satisfied. The main difficulty in proving Theorem 3 is to prove the following result.

PROPOSITION 6.1. *Let $\varepsilon > 0$. Consider a process which evolves according to the rules of model $M_1$ but starting with $[\varepsilon N]$ type 1 individuals and all other individuals having type 0. Let $g_{N,m}(\varepsilon)$ be the probability that either a type $m$ individual is born at some time or eventually all $N$ individuals have type greater than zero. Then*

$$\lim_{\varepsilon \to 0} \liminf_{N \to \infty} \varepsilon^{-1} g_{N,m}(\varepsilon) = \lim_{\varepsilon \to 0} \limsup_{N \to \infty} \varepsilon^{-1} g_{N,m}(\varepsilon) = \alpha,$$

*where $\alpha$ is as defined in (1.4).*

The first lemma will allow us to ignore overlap between type 2 families.

LEMMA 6.1. *With probability tending to one as $N \to \infty$, no type 2 individual that is born while there is an individual of type 0 in the population and another individual in the population of type 2 or higher will have a type $m$ descendant.*

PROOF. The argument is similar to the proof of Lemma 4.3. By (3.5), when we start with $[N\varepsilon]$ type 1 individuals, the total number of births and deaths of individuals of nonzero type, before the number of individuals of nonzero type reaches 0 or $N$, is at most $2\varepsilon N^2$. Since individuals give birth and die at rate 1 and mutate at rate $u_2$, the expected number of type 2 mutations is at most $\varepsilon N^2 u_2$. By (3.8), the expected amount of time during which there is an individual of type 2 or higher present in the population is at most $C\varepsilon(N^2 \log N)u_2$. Type 2 mutations happen at rate at most $Nu_2$ and produce a type $m$ descendant with probability $q_{2,m}$, so the probability that a type 2 individual born while there is another individual in the population of type 2 or higher produces a type $m$ descendant is at most $C\varepsilon N^3(\log N)u_2^2 q_{2,m}$, which is at most

(6.1) $$C(Nr_{1,m})^2(N \log N)u_2,$$

because $u_2 r_{2,m} = r_{1,m}^2$ and $q_{2,m} \sim r_{2,m}$ by Corollary 4.1. Also, we are assuming $Nr_{1,m} \to \gamma^{1/2}$, and (ii) gives $r_{1,m} \geq C u_2^{1-1/2^{m-1}}$ for some constant $C$. Therefore, $\limsup_{N \to \infty} N u_2^{1-1/2^{m-1}} < \infty$, which in combination with (iii) implies that

(6.2) $$(N \log N) u_2 \to 0.$$

It follows that the expression in (6.1) tends to zero as $N \to \infty$. □

In view of Lemma 6.1, it suffices to prove Proposition 6.1 for model $M_2$, in which no type 2 mutation can occur while there is another individual of type 2 or higher in the population. We will work with model $M_2$ for the rest of this section. As in the proof of Theorem 1, we need to deal with the correlations between individuals of type 1 and of types $j \geq 2$ caused by the fact that individuals of one positive type may replace another. To do this, we cut out the time intervals in which an individual of type 2 or higher is present in the population.

Let $X_i(t)$ be the number of type $i$ individuals at time $t$. Let

$$f(t) = \sup\left\{s : \int_0^s \mathbf{1}_{\{X_0(t) + X_1(t) = N\}} \, du = t\right\}$$

and let $Y(t) = X_1(f(t))$, so the process $(Y(t), t \geq 0)$ tracks the evolution of the number of type 1 individuals after one cuts out the times at which individuals of type $j \geq 2$ are present. Let $\beta_0 = 0$. For $i \geq 1$, let $\beta_i$ be the first time $t$ after $\beta_{i-1}$ such that $Y(t) \neq Y(t-)$ and there is no type two individual alive at time $f(t)-$, assuming such a time exists which it will a.s. as long as $Y(\beta_{i-1}) \notin \{0, N\}$. That is, the times $\beta_i$ are the times of $Y$ process jumps that happen because of a birth or death of a type one individual and do not involve the birth of a type two individual. Let $g(t) = \max\{i : \beta_i \leq t\}$, so $g(t)$ is the number of these jumps that have happened by time $t$.

We now define a discrete-time process $(Z_i)_{i=0}^\infty$, which omits the jumps in $Y$ due to time intervals being removed, but retains all of the other jumps of size 1. Let $Z_0 = [N\varepsilon]$. If $i \geq 1$, $Y(\beta_{i-1}) \notin \{0, N\}$, and $\varepsilon^3 N < Z_{i-1} < (1 - \varepsilon^2)N$, then let $Z_i = Z_{i-1} + 1$ if $Y(\beta_i) = Y(\beta_i -) + 1$, and let $Z_i = Z_{i-1} - 1$ if $Y(\beta_i) = Y(\beta_i -) - 1$. Using this induction, we can define the process $(Z_i)_{i=0}^T$, where $T = \inf\{i : Y(\beta_i) \in \{0, N\}, Z_i \leq \varepsilon^3 N, \text{ or } Z_i \geq (1-\varepsilon^2)N\}$. On the event that $\varepsilon^3 N < Z_{i-1} < (1 - \varepsilon^2)N$ and $0 < Y(\beta_i) < N$, we have $P(Z_i = Z_{i-1} + 1 | Z_0, \ldots, Z_{i-1}) = P(Z_i = Z_{i-1} - 1 | Z_0, \ldots, Z_{i-1}) = 1/2$. We then continue the process for $i > T$ by setting $Z_i$ to be $Z_{i-1} + 1$ or $Z_{i-1} - 1$ with probability $1/2$ each, independently of the population process. The process $(Z_i)_{i=0}^\infty$ is therefore a simple random walk.

Note that $T$ is smaller than the absorption time of the process $(Z_i)_{i=0}^\infty$ in $\{0, N\}$, which can be compared to the absorption time of model $M_0$ started with $[N\varepsilon]$ type 1 individuals. It therefore follows from (3.9) that $E[\beta_T] \leq N$. Thus, if $\theta > 0$, then by Markov's inequality,

$$P\left(\beta_T > \frac{N}{\theta}\right) \leq \theta. \tag{6.3}$$

Likewise, since $T$ is at most the number of births and deaths of individuals of nonzero type started from $[N\varepsilon]$ such individuals, (3.5) gives $E[T] \leq 2N^2\varepsilon \leq 2N^2$. Therefore, for $\theta > 0$,

$$P\left(T > \frac{2N^2}{\theta}\right) \leq \theta. \tag{6.4}$$

LEMMA 6.2. *For all $\delta > 0$, we have*

$$\lim_{N \to \infty} P\left(\max_{0 \leq t \leq \beta_T} |Y(t) - Z_{g(t)}| > \delta N\right) = 0.$$

PROOF. Let $\zeta_0 = 0$ and for $i \geq 1$, let $\zeta_i$ be the first time $t$ after $\zeta_{i-1}$ such that there is a type 2 individual alive at time $f(t)-$, provided such a time exists. Thus, the times $\zeta_i$ for $i \geq 1$ are the times at which the process $(Y(t), t \geq 0)$ possibly jumps because we have cut out the lifetime of a type 2 family. Every jump time of $(Y(t), t \geq 0)$ is either $\beta_i$ or $\zeta_i$ for some $i$. Since only the jumps at the times $\beta_i$ are incorporated into the process $(Z_i)_{i=1}^{\infty}$, we have

$$Y(t) - Z_{g(t)} = \sum_{i : \zeta_i \leq t} (Y(\zeta_i) - Y(\zeta_i -)). \tag{6.5}$$

We will show that the right-hand side is small because type 2 individuals are not alive in the population for a long enough time for large changes in the size of the type 1 population to happen during this time.

A type 1 individual is lost whenever a type 2 individual is born. The other changes in the number of type 1 individuals that contribute to the right-hand side of (6.5) are births and deaths that occur while there are already type 2 individuals in the population. Let $\xi_i = 1$ if the $i$th such event is a birth, and let $\xi_i = -1$ if the $i$th such event is a death. Let $J$ be the number of such events before time $f(\beta_T)$, so if $S_j = \xi_1 + \cdots + \xi_J$, then

$$|Y(t) - Z_{g(t)}| \leq |\{i : \zeta_i \leq T\}| + \max_{j \leq J} |S_j| \tag{6.6}$$

for all $t \leq \beta_T$.

The first term on the right-hand side of (6.6) is the number of type 2 mutations by time $\beta_T$, so as noted above its expected value is at most $\varepsilon N^2 u_2$. It follows from Markov's inequality and (6.2) that $P(|\{i : \zeta_i \leq T\}| > \delta N/2) \leq 4\varepsilon N^2 u_2/(\delta N) \to 0$ as $N \to \infty$.

Since $(S_j)_{j=1}^{\infty}$ is a simple random walk, by the monotone convergence theorem, the $L^2$-maximal inequality for martingales, and Wald's second equation, we have

$$E\left[\max_{j \leq J} S_j^2\right] = \lim_{n \to \infty} E\left[\max_{j \leq J \wedge n} S_j^2\right] \leq 4 \lim_{n \to \infty} E[S_{J \wedge n}^2]$$
$$= 4 \lim_{n \to \infty} E[J \wedge n] = 4E[J].$$

We have observed that the expected amount of time for which there is an individual of type 2 or greater present in the population is at most $C\varepsilon(N^2 \log N)u_2$. The rate at which type one individuals are either being born or dying is always at most $2N$, so $E[J] \leq 2C\varepsilon(N^3 \log N)u_2$. By Chebyshev's inequality and (6.2),

$$\limsup_{N \to \infty} P\left(\max_{j \leq J} |S_j| > \frac{\delta N}{2}\right) \leq \limsup_{N \to \infty} \frac{16 E[J]}{\delta^2 N^2}$$
$$\leq \limsup_{N \to \infty} \frac{32 C\varepsilon (N \log N) u_2}{\delta^2} = 0$$

and the result follows.  $\square$

LEMMA 6.3.  *For all $\delta > 0$, we have*

$$\lim_{N \to \infty} P\left(\left|\int_0^{\beta_T} (Y(t) - Z_{g(t)}) \, dt\right| > \delta N^2\right) = 0.$$

PROOF.  Let $\theta > 0$. By Lemma 6.2 and (6.3),

$$\limsup_{N \to \infty} P\left(\left|\int_0^{\beta_T} (Y(t) - Z_{g(t)}) \, dt\right| > \delta N^2\right)$$
$$\leq \limsup_{N \to \infty} \left(P\left(\beta_T > \frac{N}{\theta}\right) + P\left(\max_{0 \leq t \leq \beta_T} |Y(t) - Z_{g(t)}| > \delta \theta N\right)\right) \leq \theta.$$

Letting $\theta \to 0$ gives the result.  $\square$

LEMMA 6.4.  *For all $\delta > 0$, we have*

$$\lim_{N \to \infty} P\left(\left|\int_0^{\beta_T} Z_{g(t)} \, dt - \sum_{i=0}^{T-1} \frac{N}{2(N - Z_i)}\right| > \delta N^2\right) = 0.$$

PROOF.  For $i \leq T - 1$, let

$$D_i = \frac{N}{2(N - Z_i)} - (\beta_{i+1} - \beta_i) Z_i.$$

We need to show that

(6.7) $$\lim_{N\to\infty} P\left(\left|\sum_{i=0}^{T-1} D_i\right| > \delta N^2\right) = 0.$$

At time $t$, events that cause the number of type 1 individuals to change but do not involve the birth of a type 2 happen at rate $2Y(t)(N - Y(t))/N$. Therefore, if we define

$$\xi_i = \int_{\beta_i}^{\beta_{i+1}} \frac{2Y(t)(N - Y(t))}{N}\, dt,$$

then the random variables $\xi_i$ are independent and have the exponential distribution with mean one. Note that the process $Y$ is constant on the intervals $(\beta_i, \beta_{i+1})$ except when type 2 mutations occur. For $i \leq T - 1$, let

$$\tilde{D}_i = \frac{N}{2(N - Z_i)}(1 - \xi_i).$$

Let $\theta > 0$, so $P(T > 2N^2/\theta) \leq \theta$ by (6.4). For $0 \leq j \leq [2N^2/\theta]$, let $M_j = \sum_{i=0}^{(T-1)\wedge j} \tilde{D}_i$. Let $\mathcal{F}_j$ be the $\sigma$-field generated by $(Y(t), 0 \leq t \leq \beta_j)$. Note that $E[\tilde{D}_i|\mathcal{F}_i] = 0$, so the process $(M_j)_{j=0}^{[2N^2/\theta]}$ is a martingale. On the event that $i \leq T - 1$, we have $Z_i \leq (1 - \varepsilon^2)N$, and hence

$$\mathrm{Var}(\tilde{D}_i|\mathcal{F}_i) = \frac{N^2}{4(N - Z_i)^2} \leq \frac{1}{4\varepsilon^4}.$$

It follows from the $L^2$-maximal inequality for martingales, and orthogonality of martingale increments that

$$E\left(\max_{0\leq j\leq [2N^2/\theta]} M_j^2\right) \leq 4E[M^2_{[2N^2/\theta]}] \leq 4 \cdot \frac{2N^2}{\theta} \cdot \frac{1}{4\varepsilon^4} = \frac{2N^2}{\theta\varepsilon^4}.$$

Using Chebyshev's inequality,

$$P\left(\left|\sum_{i=0}^{T-1} \tilde{D}_i\right| > \frac{\delta N^2}{2}\right) \leq \theta + P\left(\max_{0\leq j\leq [2N^2/\theta]} |M_j| > \frac{\delta N^2}{2}\right)$$

$$\leq \theta + \frac{4}{\delta^2 N^4}\left(\frac{2N^2}{\theta\varepsilon^4}\right) = \theta + \frac{8}{\theta\delta^2\varepsilon^4 N^2}.$$

Since $\theta > 0$ was arbitrary, it follows that

(6.8) $$\lim_{N\to\infty} P\left(\left|\sum_{i=0}^{T-1} \tilde{D}_i\right| > \frac{\delta N^2}{2}\right) = 0.$$

To convert this into a bound on the $D_i$, we note that

$$|D_i - \tilde{D}_i| = \left| \frac{N}{2(N-Z_i)} \int_{\beta_i}^{\beta_{i+1}} \frac{2Y(t)(N-Y(t))}{N} \, dt - (\beta_{i+1} - \beta_i) Z_i \right|$$

$$\leq \int_{\beta_i}^{\beta_{i+1}} \left| \frac{Y(t)(N-Y(t))}{N-Z_i} - Z_i \right| dt.$$

On the event that $|Y(t) - Z_{g(t)}| \leq \gamma N$ for all $0 \leq t \leq \beta_T$, there is a constant $C_\varepsilon$ depending on $\varepsilon$ such that for all $i \leq T-1$ and $t \in [\beta_i, \beta_{i+1}]$, we have

$$\frac{Y(t)(N-Y(t))}{N-Z_i} - Z_i \leq \frac{(Z_i + \gamma N)(N - Z_i + \gamma N)}{N - Z_i} - Z_i$$

$$\leq (Z_i + \gamma N)\left(1 + \frac{\gamma}{\varepsilon^2}\right) - Z_i \leq C_\varepsilon \gamma N,$$

where in the second inequality we have used $Z_i \leq (1 - \varepsilon^2)N$. For a bound in the other direction, we note that

$$\frac{Y(t)(N-Y(t))}{N-Z_i} - Z_i \geq \frac{(Z_i - \gamma N)(N - Z_i - \gamma N)}{N - Z_i} - Z_i$$

$$\geq (Z_i - \gamma N)\left(1 - \frac{\gamma}{\varepsilon^2}\right) - Z_i \geq -C_\varepsilon \gamma N.$$

Thus, if we let $\theta > 0$ and $\gamma = \delta \theta / 2 C_\varepsilon$, then for sufficiently large $N$,

$$P\left( \left| \sum_{i=0}^{T-1} (D_i - \tilde{D}_i) \right| > \frac{\delta N^2}{2} \right)$$

$$\leq P\left( \beta_T > \frac{N}{\theta} \right) + P\left( \max_{0 \leq t \leq \beta_T} |Y(t) - Z_{g(t)}| > \gamma N \right).$$

Using (6.3), Lemma 6.2, and letting $\theta \to 0$, we get

(6.9) $$\lim_{N \to \infty} P\left( \left| \sum_{i=0}^{T-1} (D_i - \tilde{D}_i) \right| > \frac{\delta N^2}{2} \right) = 0.$$

Now (6.7) follows from (6.8) and (6.9). □

Let $D$ be the event that either $Z_T \geq (1 - \varepsilon^2)N$ or some type 2 mutation that occurs before time $f(\beta_T)$ has a type $m$ descendant.

LEMMA 6.5. *We have*

$$\lim_{N \to \infty} \left( (1 - P(D)) - E\left[ \exp\left( -r_{2,m} \sum_{i=0}^{T-1} \frac{u_2 N}{2(N - Z_i)} \right) \mathbf{1}_{\{Z_T \leq \varepsilon^3 N\}} \right] \right) = 0.$$

PROOF. If there is no type 2 individual in the population at time $t$, then the rate at which a type 2 individual is born is $u_2 X_1(t)$. Because no type 2 mutations occur while there is another type 2 individual in the population, each mutant type 2 individual independently has a type $m$ descendant with probability $q_{2,m}$. It follows that there is a mean one exponential random variable $\xi$ such that some original type two individual born before time $f(\beta_T)$ has a type $m$ descendant if and only if

$$\xi \leq \int_0^{\beta_T} Y(t) u_2 q_{2,m}\, dt. \tag{6.10}$$

Because changes in the population resulting from the birth of a type 2 individual are not recorded in the process $(Z_i)_{i=0}^{T-1}$, the random variable $\xi$ can be constructed to be independent of the process $(Z_i)_{i=0}^{T-1}$. Therefore, by conditioning on $(Z_i)_{i=0}^{T-1}$, we get

$$\begin{aligned}
P&\left(\{Z_T \leq \varepsilon^3 N\} \cap \left\{\xi > r_{2,m} \sum_{i=0}^{T-1} \frac{u_2 N}{2(N-Z_i)}\right\}\right) \\
&= E\left[\exp\left(-r_{2,m} \sum_{i=0}^{T-1} \frac{u_2 N}{2(N-Z_i)}\right) \mathbf{1}_{\{Z_T \leq \varepsilon^3 N\}}\right].
\end{aligned} \tag{6.11}$$

The event that $D$ fails to occur is the same as the event that $Z_T \leq \varepsilon^3 N$ and that (6.10) fails to occur. It follows that the difference between $P(D^c) = 1 - P(D)$ and the probability in (6.11) is at most the probability that $\xi$ is between $\int_0^{\beta_T} Y(t) u_2 q_{2,m}\, dt$ and $r_{2,m} \sum_{i=0}^{T-1} u_2 N/(2(N-Z_i))$. To bound the difference between these quantities, note that Lemmas 6.3 and 6.4 give

$$\lim_{N \to \infty} P\left(\left|\int_0^{\beta_T} u_2 Y(t)\, dt - \sum_{i=0}^{T-1} \frac{u_2 N}{2(N-Z_i)}\right| > \delta N^2 u_2\right) = 0$$

for all $\delta > 0$. Since $r_{1,m}^2 = u_2 r_{2,m}$ and $(Nr_{1,m})^2 \to \gamma$, we see that $N^2 u_2 r_{2,m}$ stays bounded as $N \to \infty$ and it follows that

$$\lim_{N \to \infty} P\left(\left|\int_0^{\beta_T} u_2 r_{2,m} Y(t)\, dt - r_{2,m} \sum_{i=0}^{T-1} \frac{u_2 N}{2(N-Z_i)}\right| > \frac{\delta}{2}\right) = 0 \tag{6.12}$$

for all $\delta > 0$. Also, $q_{2,m} \sim r_{2,m}$ by Corollary 4.1 and $P(\beta_T > N/\theta) \leq \theta$ by (6.3). Since $N^2 u_2 r_{2,m}$ stays bounded,

$$\begin{aligned}
\limsup_{N \to \infty} &P\left(\left|\int_0^{\beta_T} u_2 q_{2,m} Y(t)\, dt - \int_0^{\beta_T} u_2 r_{2,m} Y(t)\, dt\right| > \frac{\delta}{2}\right) \\
&\leq \limsup_{N \to \infty} P\left(N u_2 \beta_T |r_{2,m} - q_{2,m}| > \frac{\delta}{2}\right) = 0.
\end{aligned} \tag{6.13}$$

Combining (6.12) and (6.13) gives

$$\lim_{N \to \infty} P\left(\left|\int_0^{\beta_T} u_2 q_{2,m} Y(t)\, dt - r_{2,m} \sum_{i=0}^{T-1} \frac{u_2 N}{2(N - Z_i)}\right| > \delta\right) = 0.$$

Since

$$P\left(r_{2,m} \sum_{i=0}^{T-1} \frac{u_2 N}{2(N - Z_i)} - \delta \leq \xi \leq r_{2,m} \sum_{i=0}^{T-1} \frac{u_2 N}{2(N - Z_i)} + \delta\right) \leq 2\delta,$$

it follows that

$$\limsup_{N \to \infty} \left|(1 - P(D)) - E\left[\exp\left(-r_{2,m} \sum_{i=0}^{T-1} \frac{u_2 N}{2(N - Z_i)}\right) \mathbf{1}_{\{Z_T \leq \varepsilon^3 N\}}\right]\right| \leq 2\delta,$$

and the result follows by letting $\delta \to 0$.  □

Let $A$ be the event that either $Y(t) = N$ for some $t$, or a type $m$ individual is born at some time.

LEMMA 6.6. *There exists a constant $C$, not depending on $\varepsilon$ or $N$, such that*

$$|P(A) - P(D)| \leq C\varepsilon^2.$$

PROOF. Let $\delta > 0$, and assume that $|Y(t) - Z_{g(t)}| \leq \delta N$ for $0 \leq t \leq \beta_T$. First, suppose $D$ occurs. If a type 2 mutation that occurs before time $f(\beta_T)$ has a type $m$ descendant, then $A$ must occur. If $Z_T \geq (1 - \varepsilon^2)N$, then $Y(\beta_T) \geq (1 - \varepsilon^2 - \delta)N$, and conditional on this event the probability that $Y(t) = N$ for some $t$, in which case $A$ occurs, is at least $1 - \varepsilon^2 - \delta$. Therefore, using Lemma 6.2,

$$\limsup_{N \to \infty} P(D \cap A^c) \leq \varepsilon^2 + \delta.$$

Now, suppose $D^c$ occurs. Note that if $\delta < \varepsilon^3$ and $|Y(t) - Z_{g(t)}| \leq \delta N$ for $0 \leq t \leq \beta_T$, then we cannot have $Y(\beta_T) \in \{0, N\}$, which means we must have $Z_T \leq \varepsilon^3 N$ and, therefore, $Y(\beta_T) \leq (\varepsilon^3 + \delta)N$. Conditional on this event, the probability that $Y(t) = N$ for some $t$ is at most $\varepsilon^3 + \delta$, and the probability that one of the type one individuals at time $f(\beta_T)$ has a type $m$ descendant is at most $(\varepsilon^3 + \delta)N q_{1,m}$. From these bounds and Lemma 6.2, it follows that

$$\limsup_{N \to \infty} P(D^c \cap A) \leq (1 + \gamma^{1/2})(\varepsilon^3 + \delta).$$

The lemma follows by letting $\delta \to 0$.  □

Now let $(B_t)_{t \geq 0}$ be a Brownian motion with $B_0 = \varepsilon$. Let $U = \inf\{t : B_t = \varepsilon^3 \text{ or } B_t = 1 - \varepsilon^2\}$.

LEMMA 6.7. *We have*

$$\lim_{N \to \infty} E\left[\exp\left(-r_{2,m} \sum_{i=0}^{T-1} \frac{u_2 N}{2(N - Z_i)}\right) \mathbf{1}_{\{Z_T \leq \varepsilon^3 N\}}\right]$$
$$= E\left[\exp\left(-\frac{\gamma}{2} \int_0^U \frac{1}{1 - B_t} dt\right) \mathbf{1}_{\{B_U = \varepsilon^3\}}\right].$$

PROOF. Define a process $(W_t)_{t \geq 0}$ such that $W_t = N^{-1} Z_{[N^2 t]}$. Let $R = \inf\{t : W_t \leq \varepsilon^3 \text{ or } W_t > 1 - \varepsilon^2\}$. Note that $R = T/N^2$ and $\mathbf{1}_{\{Z_T \leq \varepsilon^3 N\}} = \mathbf{1}_{\{W_R \leq \varepsilon^3\}}$ on the event that for some $\delta < \varepsilon^3$, we have $|Y(t) - Z_{g(t)}| \leq \delta N$ for $0 \leq t \leq \beta_T$, which by Lemma 6.2 happens with probability tending to one as $N \to \infty$.

Let $\delta < \varepsilon^3$. For random variables $X_N^{(1)}$ and $X_N^{(2)}$, write $X_N^{(1)} \approx X_N^{(2)}$ if for all $\eta > 0$, there is an $N(\eta)$ such that if $N \geq N(\eta)$ then $|X_N^{(1)}/X_N^{(2)} - 1| < \eta$ on the event that $|Y(t) - Z_{g(t)}| \leq \delta N$ for $0 \leq t \leq \beta_T$. We have

$$\frac{1}{2} \int_0^R \frac{1}{1 - W_t} dt \approx \frac{1}{2} \int_0^R \frac{1}{1 - N^{-1} Z_{[N^2 t]}} dt = \frac{1}{2} \int_0^{N^2 R} \frac{1}{1 - N^{-1} Z_{[s]}} N^{-2} ds$$
$$= N^{-2} \int_0^T \frac{N}{2(N - Z_{[s]})} ds = N^{-2} \sum_{i=0}^{T-1} \frac{N}{2(N - Z_i)}.$$

Since $u_2 r_{2,m} = r_{1,m}^2$ and $(N r_{1,m})^2 \to \gamma$, we have

$$r_{2,m} \sum_{i=0}^{T-1} \frac{u_2 N}{2(N - Z_i)} \approx \gamma N^{-2} \sum_{i=0}^{T-1} \frac{N}{2(N - Z_i)} \approx \frac{\gamma}{2} \int_0^R \frac{1}{1 - W_t} dt.$$

In view of Lemma 6.2, it follows that

$$\lim_{N \to \infty} \left( E\left[\exp\left(-r_{2,m} \sum_{i=0}^{T-1} \frac{u_2 N}{2(N - Z_i)}\right) \mathbf{1}_{\{Z_T \leq \varepsilon^3 N\}}\right] \right.$$
$$\left. - E\left[\exp\left(-\frac{\gamma}{2} \int_0^R \frac{1}{1 - W_t} dt\right) \mathbf{1}_{\{W_R = \varepsilon^3\}}\right] \right) = 0.$$

Thus, it suffices to show that for all $\lambda > 0$, we have

(6.14)
$$\lim_{N \to \infty} E\left[\exp\left(-\lambda \int_0^R \frac{1}{1 - W_t} dt\right) \mathbf{1}_{\{W_R = \varepsilon^3\}}\right]$$
$$= E\left[\exp\left(-\lambda \int_0^U \frac{1}{1 - B_t} dt\right) \mathbf{1}_{\{B_U = \varepsilon^3\}}\right].$$

Since $(Z_i)_{i=0}^{\infty}$ is a simple random walk, $(W_t)_{0 \leq t \leq s}$ converges weakly as $N \to \infty$ to $(B_t)_{0 \leq t \leq s}$ for all $s > 0$. Let $D[0, s]$ be the set of real-valued functions defined on $[0, s]$ which are right continuous and have left limits. If $g:D[0,s]\to\mathbb{R}$ is bounded, and if the set of points at which it is not continuous has Wiener measure zero, then the weak convergence of $(W_t)_{0\leq t\leq s}$ to $(B_t)_{0\leq t\leq s}$ implies that $\lim_{N\to\infty} E[g((W_t)_{0\leq t\leq s})] = E[g((B_t)_{0\leq t\leq s})]$. Therefore,

$$\lim_{N\to\infty} E\left[\exp\left(-\lambda \int_0^{R\wedge s} \frac{1}{1-W_t}\,dt\right)\mathbf{1}_{\{W_{R\wedge s}=\varepsilon^3\}}\right]$$
(6.15)
$$= E\left[\exp\left(-\lambda \int_0^{U\wedge s} \frac{1}{1-B_t}\,dt\right)\mathbf{1}_{\{B_{U\wedge s}=\varepsilon^3\}}\right].$$

Note that if $\omega:[0,s]\to\mathbb{R}$ is continuous, then the function $g$ used in (6.15) is continuous at $\omega$ unless either $\inf\{t:\omega(t)=\varepsilon^3\}<\inf\{t:\omega(t)<\varepsilon^3\}$ or $\inf\{t:\omega(t)=1-\varepsilon^2\}<\inf\{t:\omega(t)>1-\varepsilon^2\}$, which would happen if $\omega$ reaches a local minimum when it first hits $\varepsilon^3$ or a local maximum when it first hits $1-\varepsilon^2$. Brownian motion paths almost surely do not have this property, so (6.15) is valid. Finally, (6.14) follows from (6.15) by letting $s\to\infty$. $\square$

Let $V = \inf\{t:B_t=0 \text{ or } B_t=1\}$.

LEMMA 6.8. *Let $I(s)=\int_0^s \frac{1}{1-B_t}\,dt$. If $\lambda>0$, there is a constant $C$ such that*

(6.16) $|E[\exp(-\lambda I(U))\mathbf{1}_{\{B_U=\varepsilon^3\}}] - E[\exp(-\lambda I(V))\mathbf{1}_{\{B_V=0\}}]| \leq C\varepsilon^2.$

PROOF. Define a process $(B'_t)_{t\geq 0}$ by $B'_t = B_{U+t}$. Let $\tau'_a = \inf\{t:B'_t=a\}$. Let $D_1$ be the event that $B_U = 1-\varepsilon^2$ and $B_V = 0$. Let $D_2$ be the event that $B_U = \varepsilon^3$ and $\tau'_{1/2} < \tau'_0$. Let $D_3$ be the event that $B_U = \varepsilon^3$ and $\tau'_0 > \varepsilon^2$. Note that on the event $(D_1 \cup D_2 \cup D_3)^c$, we have $\mathbf{1}_{\{B_U=\varepsilon^3\}} = \mathbf{1}_{\{B_V=0\}}$ and on this event we have

$$0 \leq \int_0^V \frac{1}{1-B_t}\,dt - \int_0^U \frac{1}{1-B_t}\,dt \leq 2(V-U) \leq 2\varepsilon^2.$$

It follows that the left-hand side of (6.16) is at most $P(D_1) + P(D_2) + P(D_3) + 2\lambda\varepsilon^2$.

Because Brownian motion is a martingale, we have $P(D_1) \leq P(B_V = 0|B_U = 1-\varepsilon^2) = \varepsilon^2$ and likewise $P(D_2) \leq 2\varepsilon^3$. Therefore, it remains only to show that $P(D_3) \leq C\varepsilon^2$. By the reflection principle,

$$\tfrac{1}{2} P(\tau'_0 \leq \varepsilon^2 | B_U = \varepsilon^3) = P(B'_{\varepsilon^2} \leq 0).$$

Also, $P(B'_{\varepsilon^2} > \varepsilon^3 | B_U = \varepsilon^3) = 1/2$. Therefore, $P(0 < B'_{\varepsilon^2} < \varepsilon^3 | B_U = \varepsilon^3) = [1 - P(\tau'_0 \leq \varepsilon^2 | B_U = \varepsilon^3)]/2$. It follows that

$$P(D_3) \leq P(\tau'_0 > \varepsilon^2 | B_U = \varepsilon^3)$$
$$= 2P(0 < B'_{\varepsilon^2} < \varepsilon^3 | B_U = \varepsilon^3)$$
$$\leq \frac{2\varepsilon^3}{\sqrt{2\pi\varepsilon^2}} = \varepsilon^2 \sqrt{\frac{2}{\pi}}$$

and the result follows. □

LEMMA 6.9. *Let $E_x$ denote expectation for the Brownian motion $(B_t)_{t\geq 0}$ starting from $B_0 = x$. Let*

$$u(x) = E_x\left[\exp\left(-\frac{\gamma}{2}\int_0^V \frac{1}{1-B_t}\,dt\right)\mathbf{1}_{\{B_V=0\}}\right].$$

*Then $\lim_{x\to 0} x^{-1}(1-u(x)) = \alpha$, where $\alpha$ is as defined in (1.4).*

PROOF. We choose $f$ so that $f(0) = 1$ and $f(1) = 0$. Let $g(x) = \gamma/[2(1-x)]$. Then for $0 < x < 1$, we have $u(x) = E_x[f(B_V)\exp(-\int_0^V g(B_s)\,ds)]$. Clearly $u(0) = 1$ and $u(1) = 0$. By the Feynman–Kac formula (see (6.3) on page 161 of [6]), if $v:[0,1]\to\mathbb{R}$ is a bounded continuous function such that $v(0) = 1$, $v(1) = 0$, and $\frac{1}{2}v''(x) - g(x)v(x) = 0$ for $x \in (0,1)$, then $u(x) = v(x)$ for $x \in [0,1]$. Note that (6.3) on page 161 of [6] requires $g$ to be bounded on $(0,1)$, which it is not in this example. However, the result nevertheless holds because $g$ is nonnegative and, therefore, $\exp(-\int_0^t g(B_s)\,ds)$ is always in $[0,1]$.

Multiplying by $2(1-x)$, we can write the differential equation above as $(1-x)v''(x) - \gamma v(x) = 0$. Let

(6.17) $$v(x) = c\sum_{k=1}^{\infty} \frac{\gamma^k}{k!(k-1)!}(1-x)^k,$$

where $c = 1/\sum_{k=1}^{\infty} \gamma^k/k!(k-1)!$. Note that $v(0) = 1$ and $v(1) = 0$. The series converges absolutely and uniformly on all compact subsets of $\mathbb{R}$ and can be differentiated twice term by term, so

$$(1-x)v''(x) = c\sum_{k=2}^{\infty} \frac{\gamma^k}{k!(k-1)!}k(k-1)(1-x)^{k-1}.$$

Therefore,

$$(1-x)v''(x) - \gamma v(x) = c\sum_{k=1}^{\infty}\left(\frac{\gamma^{k+1}}{k!(k-1)!}(1-x)^k - \frac{\gamma^{k+1}}{k!(k-1)!}(1-x)^k\right) = 0.$$

Thus, $v(x) = u(x)$ for $x \in [0,1]$. From our formula, it follows that

$$\lim_{x\to 0}\frac{1-u(x)}{x} = -u'(0) = c\sum_{k=1}^{\infty}\frac{\gamma^k}{(k-1)!(k-1)!} = \alpha,$$

as claimed. □

PROOF OF PROPOSITION 6.1. The only difference between $g_{N,j}(\varepsilon)$ and $P(A)$ is that the event $A$ is defined using model $M_2$, in which new type two individuals cannot be born while there is an existing individual of type 2 or higher in the population. Therefore, it follows from Lemma 4.3 that $|P(A) - g_{N,j}(\varepsilon)| \ll [N\varepsilon] r_{1,m}$ and, therefore,
$$\lim_{N \to \infty} |P(A) - g_{N,j}(\varepsilon)| = 0$$
for all $\varepsilon > 0$. By Lemmas 6.5, 6.6, 6.7 and 6.8,
$$\limsup_{N \to \infty} |P(A) - (1 - u(\varepsilon))| \leq C\varepsilon^2.$$
Combining these results and multiplying both sides by $\varepsilon^{-1}$ gives
$$\limsup_{N \to \infty} |\varepsilon^{-1} g_{N,m}(\varepsilon) - \varepsilon^{-1}(1 - u(\varepsilon))| \leq C\varepsilon.$$
Therefore, by Lemma 6.9,
$$\lim_{\varepsilon \to 0} \liminf_{N \to \infty} \varepsilon^{-1} g_{N,m}(\varepsilon) \geq \lim_{\varepsilon \to 0}(\varepsilon^{-1}(1 - u(\varepsilon)) - C\varepsilon) = \alpha,$$
$$\lim_{\varepsilon \to 0} \limsup_{N \to \infty} \varepsilon^{-1} g_{N,m}(\varepsilon) \leq \lim_{\varepsilon \to 0}(\varepsilon^{-1}(1 - u(\varepsilon)) + C\varepsilon) = \alpha$$
and the proposition follows. □

**7. Proof of Theorem 3.** With Proposition 6.1 established, the rest of the proof is routine.

LEMMA 7.1. *Consider model $M_1$, and let $q'_m$ be the probability that either a type $m$ individual is born at some time, or at some time all individuals in the population have type greater than zero. Then $\lim_{N \to \infty} N q'_m = \alpha$.*

PROOF. The probability that the number of individuals of type greater than zero reaches $[\varepsilon N]$ is $1/[\varepsilon N]$. If, at the time $T$ when the number of individuals of nonzero type reaches $[\varepsilon N]$, we change the type of all these individuals to type 1, then the probability of either getting a type $m$ individual or eventually having all $N$ individuals of type greater than zero is $g_{N,m}(\varepsilon)$. Since changing the types in this way can only reduce the probability of interest, we have
$$q'_m \geq \frac{1}{[\varepsilon N]} g_{N,m}(\varepsilon).$$

To get an upper bound, note that the probability of either having a type $m$ individual that is descended from a type 1 individual at time $T$ or having all $N$ individuals of nonzero type is at most $g_{N,m}(\varepsilon)/[\varepsilon N]$. The only possibility not accounted for is that the type $m$ individual could be descended from a type 2 individual that is born before time $T$. However, by Lemma 4.4, the proof of which is valid under our hypotheses by Corollary 4.1, the probability that a type 2 mutation that occurs while there are fewer than $\varepsilon r_{1,m}^{-1}$ individuals in the population of type 1 or higher has a type $m$ descendant is at most $C\varepsilon/N$, where we are using that $r_{1,m}$ is $O(N)$. It follows that

$$q'_m \leq \frac{1}{[\varepsilon N]} g_{N,m}(\varepsilon) + \frac{C\varepsilon}{N}.$$

The result follows from Proposition 6.1 by first letting $N \to \infty$ and then letting $\varepsilon \to 0$. $\square$

PROOF OF THEOREM 3. As in the proof of Theorem 2, call ordinary type 1 mutations type 1a, and give each individual of type greater than zero a type 1b mutation at rate $u_1$. Mutations of type 1a and 1b will both be called type 1 mutations. Let $\gamma_i$ be the time of the $i$th type 1 mutation, so the points $(\gamma_i)_{i=1}^\infty$ form a rate $Nu_1$ Poisson process on $[0, \infty)$. Define a sequence $(\zeta_i)_{i=1}^\infty$ such that $\zeta_i = 1$ if the mutation at time $\gamma_i$ is a type 1a mutation and has a type $m$ descendant in the population at some later time (which will always happen if the mutation fixates). Let $(\tilde{\zeta}_i)_{i=1}^\infty$ be a sequence of i.i.d. random variables, independent of the population process, such that $P(\tilde{\zeta}_i = 1) = q'_m$ and $P(\tilde{\zeta}_i = 0) = 1 - q'_m$ for all $i$. Let $\zeta'_i = \zeta_i$ if all individuals at time $\gamma_i-$ have type 0, and let $\zeta'_i = \tilde{\zeta}_i$ otherwise. Let $\sigma'_m = \inf\{\gamma_i : \zeta'_i = 1\}$. It is clear from the construction that $\sigma'_m$ has the exponential distribution with rate $Nu_1 q'_m$, so Lemma 7.1 gives

(7.1) $$\lim_{N \to \infty} P(u_1 \sigma'_m > t) = \exp(-\alpha t).$$

Let $\sigma_m = \inf\{\gamma_i : \zeta_i = 1\}$, which is the first time at which a type 1a mutation occurs and the individual that gets this mutation will eventually have a type $m$ descendant. We claim that $P(\sigma'_m = \sigma_m) \to 1$ as $N \to \infty$. We can only have $\sigma'_m \neq \sigma_m$ if there is a type 1 mutation at some time $\gamma_i \leq \sigma'_m$ such that not all mutations at time $\gamma_i-$ have type 0 and either $\zeta_i = 1$ or $\tilde{\zeta}_i = 1$. Note also that in this case the first such $\gamma_i$ must occur before any type 1 mutation fixates, so it suffices to consider the $\gamma_i$ that occur before any fixation. Fix $t > 0$. The expected number of type 1 mutations before time $u_1^{-1} t$ is $(Nu_1)(u_1^{-1} t) = Nt$, so by (3.8), the expected amount of time before $u_1^{-1} t$ and before any type 1 mutation fixates that there is an individual of nonzero type in the population is at most $C(N \log N)t$. Therefore, the expected number of type 1 mutations that occur before this time is at most $C(N^2 \log N) u_1 t$. If such a birth occurs at time $\gamma_i$, the probability that either $\zeta_i$ or $\tilde{\zeta}_i$ equals one is at most $2q'_m$, so

$$P(\sigma_m \neq \sigma'_m < u^{-1} t) \leq C(N^2 \log N) u_1 t q'_m \to 0,$$

where we are using that $u_1(N \log N) \to 0$ by (ii) and (6.2) and that $q'_m$ is $O(1/N)$ by Lemma 7.1. The fact that $P(\sigma'_m = \sigma_m) \to 1$ as $N \to \infty$ follows from this result and (7.1).

It remains only to show that $u_1(\tau_m - \sigma_m) \to_p 0$. When the type 1 mutation at time $\sigma_m$ does not fixate, $\tau_m - \sigma_m$ is at most the time that it takes before all descendants of the mutation die out. When this mutation fixates, then $\tau_m - \sigma_m$ includes both the time to fixation plus the time for one individual to get $m-1$ additional mutations. The probability that a given type 1 mutation takes time $\varepsilon u_1^{-1}$ to fixate or die out is at most $C u_1 \varepsilon^{-1} \log N$, so the probability that some mutation that occurs before time $u_1^{-1} t$ takes this long to fixate or die out is at most $C(N u_1)(u_1^{-1} t)(u_1 \varepsilon^{-1} \log N)$, which approaches zero as $N \to \infty$ because $u_1(N \log N) \to 0$. Finally, if a type 1 mutation fixates, then the time until a type $m$ mutation appears can be calculated using the $m-1$ case of Theorem 2 with $u_2, \ldots, u_m$ in place of $u_1, \ldots, u_{m-1}$. The hypotheses are satisfied by the arguments given in Corollary 4.1. Theorem 2 implies that the waiting time is $O(1/(Nu_2 r_{2,m}))$. However, $1/(Nu_2 r_{2,m}) \ll u_1^{-1}$ because $u_1/u_2 < b_1^{-1}$ by (ii) and $Nr_{2,m} \to \infty$ as shown in the proof of Corollary 4.1. These observations imply $u_1(\tau_m - \sigma_m) \to_p 0$, as in the proof of Theorem 2. □

**Acknowledgments.** The authors would like to thank a referee whose careful reading of the paper and many detailed comments improved the presentation of the paper.

## REFERENCES

[1] ARMITAGE, P. and DOLL, R. (1954). The age distribution of cancer and a multi-stage theory of carcinogenesis. *Brit. J. Cancer* **8** 1–12.
[2] ARRATIA, R., GOLDSTEIN, L. and GORDON, L. (1989). Two moments suffice for Poisson approximations: The Chen–Stein method. *Ann. Probab.* **17** 9–25. MR972770
[3] ATHREYA, K. B. and NEY, P. E. (1972). *Branching Processes*. Springer, New York. MR0373040
[4] BORODIN, A. N. and SALMINEN, P. (2002). *Handbook of Brownian Motion: Facts and Formulae*, 2nd ed. Birkhäuser, Boston. MR1912205
[5] CALABRESE, P., MECKLIN, J. P., JÄRVINEN, H. J., AALTONEN, L. A., TAVARÉ, S. and SHIBATA, D. (2005). Numbers of mutations to different types of colorectal cancer. *BMC Cancer* **5** 126.
[6] DURRETT, R. (1996). *Stochastic Calculus: A Practical Introduction*. CRC Press, Boca Raton, FL. MR1398879
[7] DURRETT, R. (2005). *Probability: Theory and Examples*, 3rd ed. Duxbury, Belmont, CA. MR1609153
[8] DURRETT, R. and SCHMIDT, D. (2007). Waiting for regulatory sequences to appear. *Ann. Appl. Probab.* **17** 1–32. MR2292578
[9] DURRETT, R. and SCHMIDT, D. (2008). Waiting for two mutations: With applications to regulatory sequence evolution and the limits of Darwinian evolution. *Genetics* **180** 1501–1509.
[10] ETHIER, S. N. and KURTZ, T. G. (1986). *Markov Processes: Characterization and Convergence*. Wiley, New York. MR838085
[11] EWENS, W. J. (2004). *Mathematical Population Genetics*, 2nd ed. Springer, Berlin. MR554616

[12] Iwasa, Y., Michor, F., Komarova, N. L. and Nowak, M. A. (2005). Population genetics of tumor suppressor genes. *J. Theoret. Biol.* **233** 15–23. MR2122451

[13] Iwasa, Y., Michor, F. and Nowak, M. A. (2004). Stochastic tunnels in evolutionary dynamics. *Genetics* **166** 1571–1579.

[14] Jones, S. (2008). Comparative lesion sequencing provides insights into tumor evolution. *Proc. Natl. Acad. Sci. USA* **105** 4283–4288.

[15] Knudson, A. G. (1971). Mutation and cancer: Statistical study of retinoblastoma. *Proc. Natl. Acad. Sci. USA* **68** 820–823.

[16] Kolmorogov, A. N. (1938). Zur Lösung einer biologischen Aufgabe. *Izv. NII Mat. Mekh. Tomsk. Univ.* **2** 1–6.

[17] Komarova, N. L., Sengupta, A. and Nowak, M. A. (2003). Mutation-selection networks of cancer initiation: Tumor suppressor genes and chromosomal instability. *J. Theoret. Biol.* **223** 433–450. MR2067856

[18] Luebeck, E. G. and Moolgavkar, S. H. (2002). Multistage carcinogenesis and the incidence of colorectal cancer. *Proc. Natl. Acad. Sci.* **99** 15095–15100.

[19] Moran, P. A. P. (1958). Random processes in genetics. *Proc. Cambridge Philos. Soc.* **54** 60–71. MR0127989

[20] Nowak, M. A. (2006). *Evolutionary Dynamics: Exploring the Equations of Life*. Belknap Press, Cambridge, MA. MR2252879

[21] Schinazi, R. B. (2006). The probability of treatment induced drug resistance. *Acta Biotheoretica* **54** 13–19.

[22] Schinazi, R. B. (2006). A stochastic model for cancer risk. *Genetics* **174** 545–547.

[23] Schweinsberg, J. (2008). The waiting time for *m* mutations. *Electron. J. Probab.* **13** 1442–1478. MR2438813

[24] Wodarz, D. and Komarova, N. L. (2005). *Computational Biology Of Cancer: Lecture Notes And Mathematical Modeling*. World Scientific, Singapore.

R. Durrett
Department of Mathematics
310 Malott Hall
Cornell University
Ithaca, New York 14853-4201
USA
E-mail: rtd1@cornell.edu

D. Schmidt
Institute for Mathematics
and its Applications
114 Lind Hall, 207 Church St SE
University of Minnesota
Minneapolis, Minnesota 55455
USA
E-mail: dschmidt@ima.umn.edu

J. Schweinsberg
Department of Mathematics, 0112
University of California, San Diego
9500 Gilman Drive
La Jolla, California 92093-0112
USA
E-mail: jschwein@math.ucsd.edu